\documentclass[11pt,letterpaper]{amsart}
\usepackage[margin=1in]{geometry}

\usepackage{amsmath,amssymb,amsthm,mathrsfs,bbm,comment,units,enumitem,xcolor,graphics}
\usepackage[colorlinks, linktocpage]{hyperref}
\usepackage{pgf,tikz}
\usepackage{tikz-3dplot}
\usepackage{multicol}
\usepackage{longtable}
\usetikzlibrary{positioning}
\usetikzlibrary{decorations.pathreplacing}
\usetikzlibrary{arrows,shapes}
\usetikzlibrary{arrows.meta,bending,automata}

\renewcommand{\leq}{\leqslant}
\renewcommand{\geq}{\geqslant}

\newcommand{\One}{\text{I}}
\newcommand{\Two}{\text{II}}

\newcommand{\Tonexi}{\overset{\leftarrow}{T}_{\One, \xi}}
\newcommand{\Ttwoxi}{\overset{\rightarrow}{T}_{\Two, \xi}}
\newcommand{\Sonexi}{\overset{\leftarrow}{S}_{\One, \xi}}
\newcommand{\Stwoxi}{\overset{\rightarrow}{S}_{\Two, \xi}}

\theoremstyle{definition}
\newtheorem{theorem}{Theorem}
\newtheorem{definition}{Definition}
\newtheorem{lemma}[theorem]{Lemma}
\newtheorem{corollary}[theorem]{Corollary}
\newtheorem{proposition}[theorem]{Proposition}
\newtheorem{example}{Example}

\newtheorem{remark}{Remark}

\newtheorem*{notn}{Notation}

\title{Selection Games on Hyperspaces}
\author{Christopher Caruvana}
\address{School of Sciences\\
Indiana University Kokomo\\
2300 S. Washington Street, Kokomo, IN 46902 USA}
\email{chcaru@iu.edu}
\urladdr{https://chcaru.pages.iu.edu/}
\author{Jared Holshouser}
\address{Department of Mathematics\\
Norwich University\\
158 Harmon Drive, Northfield, VT 05663 USA}
\email{JHolshou@norwich.edu}
\urladdr{https://jaredholshouser.github.io/}
\date{\today}

\begin{document}

\maketitle
\begin{abstract}
    In this paper we connect selection principles on a topological space to corresponding selection principles on one of its hyperspaces.
    We unify techniques and generalize theorems from the known results about selection principles for common hyperspace constructions.
    This includes results of Lj.D.R. Ko{\v{c}}inac, Z. Li, and others.
    We use selection games to generalize selection principles and we work with strategies of various strengths for these games.
    The selection games we work with are primarily abstract versions of the selection principles of Rothberger, Menger, and Hurewicz type, as well as games of countable fan tightness and selective separability.
    The hyperspace constructions that we work with are the Vietoris and Fell topologies, both upper and full, generated by ideals of closed sets.
    Using a new technique we are able to extend straightforward connections between topological constructs to connections between selection games related to those constructs.
    This extension process works regardless of the length of the game, the kind of selection being performed, or the strength of the strategy being considered.
\end{abstract}

\section{Introduction} \label{sec:Intro}

In \cite[pp. 10]{Scheepers2003} Scheepers proved what he referred to as monotonicity laws for selection principles.
The monotonicity laws allow one to draw implications between selection principles from the subset relations between the collections defining the selection principle (for example: if \(\mathcal A \subseteq \mathcal B\), and \(\mathcal C\) is some other collection, then \(S_1(\mathcal B,\mathcal C) \implies S_1(\mathcal A,\mathcal C)\)).
In a prior paper \cite{CCJH2020}, the current authors created a new kind of monotonicity law. Instead of requiring direct inclusion, this new law uses translation functions which can be lined up to play games.
Throughout that paper the translation functions are used to concisely prove some general theorems regarding the connections between selection games on a space and one of its corresponding spaces of continuous functions.
Theorem 40 of \cite{CaruvanaHolshouser} identifies a relationship between a class of selection games on a space and the space of its compact subsets endowed with the Vietoris topology.
In this paper we seek to elaborate on that relationship using the translation theorems from \cite{CCJH2020}.

The three main results from this paper are vaguely described below. We direct the reader to Section \ref{section:Definitions} for notation and definitions.

\begin{theorem}
    Suppose that \(\mathcal{A, B, C, D}\) are collections so that \(\bigcup \mathcal C \subseteq \bigcup \mathcal A\) and \(\bigcup \mathcal D \subseteq \bigcup \mathcal B\).
    Additionally, suppose that there exists a bijection \(\beta : \bigcup \mathcal A \to \bigcup \mathcal B\) with the following features:
    \begin{itemize}
        \item \(A \in \mathcal A \iff \beta[A] \in \mathcal B\), and
        \item \(C \in \mathcal C \iff \beta[C] \in \mathcal D\).
    \end{itemize}
    Then the selection game defined by \(\mathcal A\) and \(\mathcal C\), and the game defined by \(\mathcal B\) and \(\mathcal D\) are equivalent.
    This equivalence holds whether the games use single selections or finite selections.
\end{theorem}

\begin{theorem}
    Fix a topological space \(X\).
    Then
    \begin{enumerate}[label=(\roman*)]
        \item The generalized Rothberger/Menger game on \(X\) is equivalent to the  Selective Separability game played on the generalized upper Fell topology,
        \item
        The generalized local Rothberger/Menger game on \(X\) is equivalent to the  Countable Fan Tightness game played on the generalized upper Fell topology, and
        \item
        The generalized Hurewicz game on \(X\) is equivalent to a modified Selective Separability game played on the generalized upper Fell topology.
    \end{enumerate}
\end{theorem}

\begin{theorem}
    The previous theorem holds between \(X\) and the generalized full Fell topology if the covers of \(X\) are replaced by \(\mathcal A_F\)-covers.
\end{theorem}

Significant prior work on selection principles for the upper Fell topology and upper Vietoris topology was done by Di Maio, Ko{\v{c}}inac, and Meccariello in 2015 \cite{KocinacEtAl2005}.
In 2016, Li developed the right definitions to work on selection principles for the full Fell topology and full Vietoris topology \cite{Li2016}.
Aside from generalizing the results of \cite{KocinacEtAl2005} and \cite{Li2016} for selection principles to selection games, we also provide a unifying infrastructure which highlights the similarities, not only in the results, but in the arguments employed.

In Section \ref{section:Definitions} we provide definitions for most of the topological concepts we will be discussing. We also summarize some of the relationships between those concepts.
In Section \ref{section:Translations}, we elaborate on the strategy translation results of \cite{CCJH2020}. The section begins with the most general and technical version of the translation theorem. We then use that to generate less general but more easily applied versions of the translation theorem.

In Section \ref{section:Covers} we establish various equivalences involving different kinds of open sets and other classes of sets.
Finally, in Section \ref{section:Applications}, we tie the equivalences of Section \ref{section:Covers} to selection principles and selection games using the results of Section \ref{section:Translations}.
As an ending note, in \ref{subsection:BarIlan}, we look at some forms of tightness and determine that the results of Section \ref{section:Translations} also link to certain selection principles.

For general background and notation involving selection principles, see \cite{SakaiScheepers} and \cite{KocinacSelectedResults}.

Throughout, all ground spaces \(X\) are assumed to be Hausdorff (and non-compact).

\section{Definitions} \label{section:Definitions}

\begin{definition}
    For a topological space \(X\), we let
    \begin{itemize}
        \item
        \(\mathscr T_X\) denote the collection of all proper, non-empty open subsets of \(X\),
        \item
        \([X]^{<\omega}\) denote the collection of all non-empty finite subsets of \(X\),
        \item
        \(\mathbb K(X)\) denote the collection of all proper, non-empty compact subsets of \(X\), and
        \item
        \(\mathbb F(X)\) denote the collection of all proper, non-empty closed subsets of \(X\).
    \end{itemize}
\end{definition}
Notice that the map \(U \mapsto X \setminus U\), \(\mathscr T_X \to \mathbb F(X)\), is a bijection.
We will be using this often.

\begin{definition}
    Given a topological space \(X\), by an \emph{ideal of closed sets} we mean a collection \(\mathcal A\) consisting of proper closed sets such that
\begin{itemize}
    \item
    for all \(A , B \in \mathcal A\), if \(A \cup B \neq X\), then \(A \cup B \in \mathcal A\) and
    \item
    if \(A \in \mathcal A\) and \(B \subseteq A\) is closed, then \(B \in \mathcal A\).
\end{itemize}
\end{definition}
Examples of ideals of closed sets include
\begin{itemize}
    \item
    the collection of all finite subsets of a space,
    \item
    the collection of all compact subsets of a non-compact space,
    \item
    the collection of closed and real-bounded sets of a space,
    \item
    the collection of all closed subsets with empty interior inside a Baire space,
    \item
    the collection of all closed subsets with measure zero in a measure space, and
    \item
    the collection of all proper closed subsets of a space.
\end{itemize}

\subsection{Types of Covers}

Generally, for an open cover \(\mathscr U\) of a topological space \(X\), we say that \(\mathscr U\) is \emph{non-trivial} provided that \(X \not\in \mathscr U\).
We let \(\mathcal O_X\) denote the collection of all non-trivial open covers of \(X\).
If \(Y \subseteq X\), we say that a collection of open sets \(\mathscr U\) is a \emph{non-trivial open cover of \(Y\)} provided that \(Y \subseteq \bigcup \mathscr U\) and, for each \(U \in \mathscr U\), \(Y \not\subseteq U\).

The first general class of covers we define here generalized the idea of an \(\omega\)-cover or a \(k\)-cover.

\begin{definition}
    For a space \(X\) and an ideal of closed sets \(\mathcal A\), a non-trivial open cover \(\mathscr U\) is an \emph{\(\mathcal A\)-cover} if, for every \(A \in \mathcal A\), there exists \(U \in \mathscr U\) so that \(A \subseteq U\).
    We let \(\mathcal O(X, \mathcal A)\) denote the collection of all \(\mathcal A\)-covers of \(X\).
    We also define \(\Gamma(X, \mathcal A)\) to be all infinite non-trivial open covers \(\mathscr U\) so that for every \(A \in \mathcal{A}\), \(\{ U \in \mathscr U : A \not\subseteq U\}\) is finite.
\end{definition}

\begin{remark}
    Note that
    \begin{itemize}
        \item
		\(\mathcal O(X,[X]^{<\omega})\) is the collection of all \(\omega\)-covers of \(X\).
		\item
		\(\mathcal O(X, \mathbb K(X))\) is the collection of all \(k\)-covers of \(X\).
		\item
		\(\Gamma(X, \mathbb K(X))\) is the collection of all \(\gamma_k\)-covers of \(X\).
    \end{itemize}
\end{remark}

A kind of localized cover construction occurs naturally when studying the upper Vietoris and Fell topologies.
We generalize that in the following definition.

\begin{definition} \label{def:LocalizedCovers}
    For a space \(X\), \(Y \subseteq X\), and an ideal of closed sets \(\mathcal A\), a non-trivial open cover \(\mathscr U\) of \(Y\) is an \emph{\(\mathcal A\)-cover of Y} if, for every \(A \in \mathcal A\) with \(A \subseteq Y\), there exists \(U \in \mathscr U\) so that \(A \subseteq U\).
    We let \(\mathcal O(X, Y, \mathcal A)\) denote the collection of all \(\mathcal A\)-covers of \(Y\).
\end{definition}

\begin{remark}
    Note that \(\mathcal O(\mathbb R, (-1,1), \mathbb F(\mathbb R))\) is not the same as \(\mathcal O((-1,1), \mathbb F((-1,1)))\).
    So in general, the localized covers in Definition \ref{def:LocalizedCovers} are not the same thing as the relatively open covers.
\end{remark}

\begin{notn}
	Suppose \(\mathcal A\) is a collection of subsets of \(X\). We define
	\begin{itemize}
	    \item
	    \(c.\mathcal {A}\) to be the set \(\{ X \setminus A : A \in \mathcal A\}\).
	    \item
	    \(\neg \mathcal A\) to be the complement of \(\mathcal A\) relative to the power set of \(X\).
	\end{itemize}
\end{notn}

The notion of a groupable cover is used to provide alternate characterizations of the Gerlits-Nagy \((\ast)\) property, introduced in \cite{GerlitsNagy}, and Hurewicz property, introduced in \cite{Hurewicz}, as seen in \cite{KocinacScheepers7} and \cite{KocinacSelectedResults}.
The idea of groupable and weakly groupable are added to the concept of an \(\mathcal A\)-cover in the following definitions.
Normally, groupable covers are assumed to be countable.
Here we allow the covers to be uncountable.

\begin{definition}
    For a space \(X\) and an ideal of closed sets \(\mathcal A\), a non-trivial open cover \(\mathscr U\) is an \emph{\(\omega\)-groupable \(\mathcal A\)-cover} if there is a partition \(\varphi:\mathscr U \to |\mathscr U|\) so that
    \begin{itemize}
        \item \(\varphi^{-1}(\xi)\) is finite for all \(\xi < |\mathscr U|\), and
        \item for every \(A \in \mathcal A\), there is a \(\xi_0 < |\mathscr U|\) so that for all \(\xi \geq \xi_0\), there is a \(U_\xi \in \varphi^{-1}(\xi)\) so that \(A \subseteq U_\xi\).
    \end{itemize}
    Note that if \(\mathscr U\) is countable, then \(\omega\)-groupable corresponds to the notion of groupability as discussed in \cite{KocinacEtAl2005}.
    We include the \(\omega\) modifier to indicate that the partition consists of finite sets.
    We will observe this convention in forthcoming definitions as well.

    We let \(\mathcal O^{gp}(X, \mathcal A)\) denote the collection of all \(\omega\)-groupable \(\mathcal A\)-covers of \(X\).
\end{definition}

\begin{definition}
    For a space \(X\) and an ideal of closed sets \(\mathcal A\), a non-trivial open cover \(\mathscr U\) is a \emph{weakly \(\omega\)-groupable \(\mathcal A\)-cover} if there is a partition \(\varphi:\mathscr U \to |\mathscr U|\) so that
    \begin{itemize}
        \item \(\varphi^{-1}(\xi)\) is finite for all \(\xi < |\mathscr U|\), and
        \item for every \(A \in \mathcal A\), there is a \(\xi < |\mathscr U|\) so that \(A \subseteq \bigcup \varphi^{-1}(\xi)\).
    \end{itemize}
    Note that if \(\mathscr U\) is countable, then being weakly \(\omega\)-groupable corresponds to the notion of being weakly groupable as discussed in \cite{BabinkostovaKocinacScheepers}.

    We let \(\mathcal O^{wgp}(X, \mathcal A)\) denote the collection of all weakly \(\omega\)-groupable \(\mathcal A\)-covers of \(X\).
\end{definition}

In \cite{Li2016}, Li defined cover types \(c_V\) and \(k_F\) which captured features of the full Vietoris and Fell topologies, respectively.
Let \(\mathcal K_F(X)\) denote the \(k_F\)-covers and \(\mathcal C_V(X)\) denote the \(c_V\)-covers of \(X\).
In the following definition we establish a cover type which generalizes both \(\mathcal K_F(X)\) and \(\mathcal C_V(X)\).
\begin{definition}
    For a space \(X\) and an ideal of closed sets \(\mathcal A\), a non-trivial open cover \(\mathscr U\) is an \emph{\(\mathcal A_F\)-cover} iff for all \(A \in \mathcal A\) and all finite sequences of open sets \(V_1,\ldots,V_n \subseteq X\) with the property that \((X \setminus A) \cap V_j \neq \emptyset\) for \(1 \leq j \leq n\), there are an open set \(U \in \mathscr U\) and a finite set \(F \subseteq X\) so that
    \begin{itemize}
        \item \(A \subseteq U\),
        \item \(F \cap V_j \neq \emptyset\) for \(1 \leq j \leq n\), and
        \item \(F \cap U = \emptyset\).
    \end{itemize}
    Let \(\mathcal O_F(X,\mathcal A)\) denote the \(\mathcal A_F\)-covers of \(X\).
\end{definition}
Note that, for any collection \(\mathcal A\), \(\mathcal O_F(X,\mathcal A) \subseteq \mathcal O(X,\mathcal A)\).
We have chosen the subscript \(F\) here to reflect the added condition regarding finite sets and to remind the reader that these covers are used to study the Fell topology.

\begin{proposition}\label{prop:DefSame}
    Fix a space \(X\). Then \(\mathcal K_F(X) = \mathcal O_F(X,\mathbb K(X))\) and \(\mathcal C_V(X) = \mathcal O_F(X,\mathbb F(X))\).
\end{proposition}
\begin{proof}
    The fact that \(\mathcal K_F(X) = \mathcal O_F(X,\mathbb K(X))\) is immediate from Definition 2.1 of \cite{Li2016}.
    We show that \(\mathcal C_V(X) = \mathcal O_F(X,\mathbb F(X))\).

    Suppose \(\mathscr U \in \mathcal C_V(X)\).
    Let \(A \subseteq X\) be a closed set and \(V_1,\ldots,V_n \subseteq X\) be open sets with the property that \((X \setminus A) \cap V_j \neq \emptyset\) for \(1 \leq j \leq n\).
    Set \(W_j = (X \setminus A) \cap V_j\).
    Then \(W_1,\ldots,W_n\) are non-empty open sets, so we can find a finite set \(F \subseteq X\) and \(U \in \mathscr U\) so that \(F \cap W_j \neq \emptyset\) for \(1 \leq j \leq n\), \(\bigcap_{j=1}^n (X \setminus W_j) \subseteq U\), and \(F \cap U = \emptyset\).

    Since \(\bigcup_{j=1}^n W_j \subseteq (X \setminus A)\),
    \[
    A \subseteq X \setminus \left(\bigcup_{j=1}^n W_j \right) = \bigcap_{j=1}^n (X \setminus W_j) \subseteq U.
    \]
    This shows that \(\mathscr U \in \mathcal O_F(X,\mathbb F(X))\).

    Now suppose that \(\mathscr U \in \mathcal O_F(X,\mathbb F(X))\).
    Let \(V_1,\ldots,V_n \subseteq X\) be open sets.
    Set \(A = \bigcap_{j=1}^n (X \setminus V_j)\).
    Then we can find \(U \in \mathscr U\) and \(F \subseteq X\) finite so that \(A \subseteq U\), \(F \cap V_j \neq \emptyset\) for \(1 \leq j \leq n\),  and \(F \cap U = \emptyset\).
    Thus \(\bigcap_{j=1}^n (X \setminus V_j) \subseteq U\) and \(\mathscr U\) is a \(c_V\)-cover.
\end{proof}

Recall that all finite sets are compact and that therefore all \(k\)-covers are \(\omega\)-covers.
This is part of a more general trend that we formalize below.

\begin{definition}
    Recall that a collection \(\mathcal B\) \emph{refines} \(\mathcal A\), denoted \(\mathcal B \prec \mathcal A\), if, for every \(B \in \mathcal B\), there exists \(A \in \mathcal A\) so that \(B \subseteq A\).
\end{definition}

\begin{proposition} \label{prop:RefinementEquivalences}
    For ideals of closed sets \(\mathcal A\) and \(\mathcal B\), the following are equivalent.
    \begin{enumerate}[label=(\roman*)]
        \item \label{prop:refines}
        \(\mathcal B \prec \mathcal A\)
        \item \label{prop:refineSubset}
        \(\mathcal B \subseteq \mathcal A\)
        \item \label{prop:normalContainment}
        \(\mathcal O(X,\mathcal A) \subseteq \mathcal O(X,\mathcal B)\)
        \item \label{prop:FContainment}
        \(\mathcal O_F(X, \mathcal A) \subseteq \mathcal O_F(X,\mathcal B)\).
    \end{enumerate}
\end{proposition}
\begin{proof}
    The equivalence of \ref{prop:refines} and \ref{prop:normalContainment} was demonstrated in Lemma 6 of \cite{CCJH2020}.

    We prove the equivalence of \ref{prop:refines} and \ref{prop:refineSubset}.
    It is immediate that that \(\mathcal B \subseteq \mathcal A\) implies \(\mathcal B \prec \mathcal A\).
    Now, if \(\mathcal B \subseteq \mathcal A\) and \(B \in \mathcal B\), let \(A \in\mathcal A\) be so that \(B \subseteq A\).
    Then, since \(B\) is closed and \(\mathcal A\) is an ideal of closed sets, \(B \in \mathcal A\).

    We now prove the equivalence of \ref{prop:refineSubset} and \ref{prop:FContainment}.
    Suppose \(\mathcal B \subseteq \mathcal A\) and let \(\mathscr U \in \mathcal O_F(X , \mathcal A)\).
    Let \(B \in \mathcal B\) and \(V_1 , V_2, \ldots , V_n\) be open so that \(V_j \setminus B \neq \emptyset\) for each \(j\).
    Notice that \(B \in \mathcal A\).
    Since \(\mathscr U \in \mathcal O_F(X , \mathcal A)\) and \(B \in \mathcal A\), there exists a \(U \in \mathscr U\) with the desired properties.

    By way of the contrapositive, suppose \(\mathcal B \not\subseteq \mathcal A\).
    Let \(B \in \mathcal B\)  be so that \(B \notin \mathcal A\).
    Notice that since \(\mathcal A\) is an ideal, this means that \(B \not\subseteq A\) for all \(A \in \mathcal A\).
    Choose \(x_A \in B \setminus A\) for all \(A \in \mathcal A\) and set \(U_A = X \setminus \{x_A\}\).
    Consider
    \[
    \mathscr U = \{ U \in \mathscr T_X : (\exists A \in \mathcal A)[U \subseteq U_A] \}.
    \]
    To see that \(\mathscr U \not\in \mathcal O_F(X, \mathcal B)\), we note that \(\mathscr U \not\in \mathcal O(X,\mathcal B)\).
    To see that \(\mathscr U \in \mathcal O_F(X, \mathcal A)\), let \(A \in \mathcal A\) and \(V_1, V_2 , \ldots, V_n\) be open sets with \(V_j \setminus A \neq\emptyset\) for each \(j\).
    Let \(x_j \in V_j \setminus A\) for each \(j = 1,2, \ldots , n\).
    Then consider \(F = \{ x_A , x_1 , x_2 , \ldots , x_n \}\) and \(U = X \setminus F\).
    Observe that \(A \subseteq U\), \(F \cap V_j \neq \emptyset\) for each \(j\), and that \(F \cap U = \emptyset\).
    Since \(U \in \mathscr U\), we see that \(\mathscr U \in \mathcal O_F(X,\mathcal A)\).
\end{proof}

Just as localized \(\mathcal A\)-covers show up when studying the upper hyperspace topologies, localized versions of the \(\mathcal A_F\)-covers naturally occur when studying the full hyperspace topologies. Li defined localized versions of \(\mathcal K_F\) and \(\mathcal C_V\) in Definitions 2.2 and 2.4 of \cite{Li2016}.
Let \(\mathcal K_F(X,Y)\) denote the \(k_F\)-covers of \(Y\) and \(\mathcal C_V(X,Y)\) denote the \(c_V\)-covers of \(Y\).
The following definition generalizes these two notions.

\begin{definition}
    For a space \(X\), non-empty \(Y \subseteq X\) and an ideal of closed sets \(\mathcal A\), a non-trivial open cover \(\mathscr U\) of \(Y\) is an \emph{\(\mathcal A_F\)-cover of Y} iff for all \(A \in \mathcal A\) with \(A \subseteq Y\) and all finite sequences of open sets \(V_1,\ldots,V_n \subseteq X\) with the property that \((X \setminus Y) \cap V_j \neq \emptyset\) for \(1 \leq j \leq n\), there are an open set \(U \in \mathscr U\) and a finite set \(F \subseteq X\) so that
    \begin{itemize}
        \item \(A \subseteq U\),
        \item \(F \cap V_j \neq \emptyset\) for \(1 \leq j \leq n\), and
        \item \(F \cap U = \emptyset\).
    \end{itemize}
    Let \(\mathcal O_F(X,Y,\mathcal A)\) denote the \(\mathcal A_F\)-covers of \(Y\).
\end{definition}

\begin{proposition}
    Fix a space \(X\), non-empty \(Y \subseteq X\), and a non-trivial open cover \(\mathscr U\) of \(X\). Then \(\mathcal K_F(X,Y) = \mathcal O_F(X,Y,\mathbb K(X))\) and \(\mathcal C_V(X,Y) = \mathcal O_F(X,Y,\mathbb F(X))\).
\end{proposition}
\begin{proof}
    The fact that \(\mathcal K_F(X,Y) = \mathcal O_F(X,Y,\mathbb K(X))\) follows immediately from Definition 2.2 of \cite{Li2016}.
    So we address \(\mathcal C_V(X,Y) = \mathcal O_F(X,Y,\mathbb F(X))\).

    Suppose \(\mathscr U \in \mathcal C_V(X,Y)\).
    Let \(A\subseteq Y\) be closed and \(V_1 , V_2 , \ldots, V_n\) be open with the property that \(V_j \setminus Y \neq \emptyset\) for each \(j\).
    Set \(W_0 = X \setminus A\) and, for each \(j = 1,2,\ldots, n\), \(W_j = V_j \setminus A\).
    Notice that each \(W_j\) is a non-empty open set.
    Moreover, \(W_j \setminus Y \neq \emptyset\) for each \(j\).
    Notice that
    \[
        \bigcap_{j=0}^n (X \setminus W_j) \subseteq X \setminus W_0 = A \subseteq Y
    \]
    so we can find \(U \in \mathscr U\) and a finite set \(F \subseteq X\) so that \(\bigcap_{j=0}^n (X \setminus W_j) \subseteq U\), \(F \cap W_j \neq\emptyset\) for each \(j\), and \(F \cap U = \emptyset\).
    As \(W_j \subseteq X \setminus A\), we see that \(A \subseteq X \setminus W_j\) for each \(j\).
    Hence, \(A \subseteq \bigcap_{j=0}^n (X \setminus W_j) \subseteq U\).
    That is, \(\mathscr U \in \mathcal O_F(X,Y,\mathbb F(X))\).

    Now, suppose \(\mathscr U \in \mathcal O_F(X,Y,\mathbb F(X))\).
    Let \(V_1 , V_2 , \ldots , V_n\) be open sets with \(\bigcap_{j=1}^n (X \setminus V_j) \subseteq Y\) and \(V_j \setminus Y \neq \emptyset\) for each \(j\).
    Let \(A = \bigcap_{j=1}^n (X \setminus V_j)\) and notice that \(A\) is closed.
    Then we can find \(U \in \mathscr U\) and finite \(F \subseteq X\) so that \(A \subseteq U\), \(F \cap V_j \neq \emptyset\) for each \(j\), and \(F \cap U\).
    This establishes that \(\mathscr U \in \mathcal C_V(X,Y)\).
\end{proof}

Following Definitions 5.1 and 5.3 of \cite{Li2016}, where the classes \(\mathcal K_F^{gp}\) of \(k_F\)-groupable covers and \(\mathcal C_V^{gp}\) of \(c_V\)-groupable covers are introduced, we define the general notion of an \(F\)-groupable \(\mathcal A_F\)-cover.

\begin{definition}
    For a space \(X\), an ideal \(\mathcal A\) of closed sets, and a non-trivial open cover \(\mathscr U\) of \(X\), we say \(\mathscr U\) is \emph{\(\omega\)-\(F\)-groupable} if there is a function \(\varphi:\mathscr U \to |\mathscr U|\) so that
    \begin{itemize}
        \item \(\varphi^{-1}(\xi)\) is finite for all \(\xi < |\mathscr U|\) and
        \item for all \(A \in \mathcal A\) and all finite sequences of open sets \(V_1,\ldots,V_m \subseteq X\) with the property that \((X \setminus A) \cap V_j \neq \emptyset\) for \(1 \leq j \leq m\), there is a \(\xi_0 < |\mathscr U|\) so that for all \(\xi \geq \xi_0\), there is an open set \(U_\xi \in \varphi^{-1}(\xi)\) and a finite set \(F_\xi \subseteq X\) so that
        \begin{itemize}
            \item \(A \subseteq U_\xi\),
            \item \(F_\xi \cap V_j \neq \emptyset\) for \(1 \leq j \leq m\), and
            \item \(F_\xi \cap U_\xi = \emptyset\).
        \end{itemize}
    \end{itemize}
    Let \(\mathcal O_F^{gp}(X,\mathcal A)\) denote the \(\omega\)-\(F\)-groupable \(\mathcal A_F\)-covers of \(X\).
\end{definition}

\begin{proposition}
    Fix a space \(X\). Generally, \(\mathcal K_F^{gp}(X) \subseteq \mathcal O_F^{gp}(X,\mathbb K(X))\) and \(\mathcal C_V^{gp}(X) \subseteq \mathcal O_F^{gp}(X,\mathbb F(X))\).
    In the case that we have a countable family \(\mathscr U\) of open sets, \(\mathscr U \in \mathcal O_F^{gp}(X,\mathbb K(X)) \implies \mathscr U \in \mathcal K_F^{gp}(X)\) and \(\mathscr U \in \mathcal O_F^{gp}(X,\mathbb F(X)) \implies \mathscr U \in \mathcal C_V^{gp}(X)\).
\end{proposition}
\begin{proof}
    The fact that \(\mathcal K_F^{gp}(X) \subseteq \mathcal O_F^{gp}(X,\mathbb K(X))\) follows from the observation that Definition 5.1 of \cite{Li2016} necessitates that \(\mathscr U \in \mathcal K_F^{gp}(X)\) be countable.
    To see that \(\mathcal C_V^{gp}(X) \subseteq \mathcal O_F^{gp}(X,\mathbb F(X))\), let \(\mathscr U \in \mathcal C_V^{gp}(X)\).
    By Definition 5.3 of \cite{Li2016}, we have \(\varphi : \mathscr U \to \omega\) so that \(\varphi^{-1}(n)\) is finite for each \(n\) and some other properties hold.
    Let \(A \subseteq X\) be closed and \(V_1 , V_2 , \ldots, V_m\) be open so that \(V_j \setminus A \neq \emptyset\) for all \(j\).
    Set \(W_j = V_j \setminus A\) and notice that each \(W_j\) is a non-empty open set.
    Then we can find \(n_0 \in \omega\) so that, for each \(n \geq n_0\), there exist \(U_n \in \varphi^{-1}(n)\) and \(F_n \in [X]^{<\omega}\) so that \(\bigcap_{j=1}^m (X \setminus W_j) \subseteq U_n\), \(F_j \cap W_j \neq\emptyset\) for each \(j\), and \(F_n \cap U_n = \emptyset\).
    Since \(\bigcup_{j=1}^m W_j \subseteq X \setminus A\), we see that \(A \subseteq \bigcap_{j=1}^m (X \setminus W_j) \subseteq U_n\).
    This demonstrates that \(\mathscr U \in \mathcal O_F^{gp}(X,\mathbb F(X))\).

    If \(\mathscr U \in \mathcal O_F^{gp}(X,\mathbb K(X))\) is countable, then we see that \(\mathscr U \in \mathcal K_F^{gp}(X)\) by \cite[Definition 5.1]{Li2016}.
    So suppose \(\mathscr U \in \mathcal O_F^{gp}(X,\mathbb F(X))\) is a countable.
    Let \(\varphi : \mathscr U \to \omega\) satisfy the membership criteria for \(\mathcal O_F^{gp}(X,\mathbb F(X))\).
    Notice that \(\{\varphi^{-1}(n):n\in\omega\}\) represents \(\mathscr U\) as a union of infinitely many finite sets.
    Let \(V_1, V_2, \ldots, V_m\) be non-empty open sets and notice that \(A = \bigcap_{j=1}^m (X \setminus V_j)\) is a proper closed subset of \(X\).
    Moreover, \((X \setminus A) \cap V_j \neq \emptyset\) for each \(j\).
    Then we obtain \(n_0 \in \omega\) so that, for every \(n \geq n_0\), there exists \(U \in \varphi^{-1}(n)\) and \(F \in [X]^{<\omega}\) so that
    \[
        \bigcap_{j=1}^m (X \setminus V_j) = A \subseteq U,
    \]
    \(F \cap V_j \neq\emptyset\) for each \(j\), and \(F \cap U = \emptyset\).
    That is, \(\mathscr U \in \mathcal C_V^{gp}(X)\).
\end{proof}

\begin{definition}
    For a space \(X\), an ideal \(\mathcal A\) of closed sets, and a non-trivial open cover \(\mathscr U\) of \(X\), we say \(\mathscr U\) is \emph{weakly \(\omega\)-\(F\)-groupable} if there is a function \(\varphi:\mathscr U \to |\mathscr U|\) so that
    \begin{itemize}
        \item \(\varphi^{-1}(\xi)\) is finite for all \(\xi < |\mathscr U|\) and
        \item for all \(A \in \mathcal A\) and all finite sequences of open sets \(V_1,\ldots,V_m \subseteq X\) with the property that \((X \setminus A) \cap V_j \neq \emptyset\) for \(1 \leq j \leq m\), there is a \(\xi < |\mathscr U|\) and a finite set \(F \subseteq X\) so that
        \begin{itemize}
            \item \(A \subseteq \bigcup \varphi^{-1}(\xi)\),
            \item \(F \cap V_j \neq \emptyset\) for \(1 \leq j \leq m\), and
            \item \(F \cap \bigcup \varphi^{-1}(\xi) = \emptyset\).
        \end{itemize}
    \end{itemize}
    Let \(\mathcal O_F^{wgp}(X,\mathcal A)\) denote the weakly \(\omega\)-\(F\)-groupable \(\mathcal A_F\)-covers of \(X\).
\end{definition}

\begin{proposition}
    Fix a space \(X\). Then \(\mathcal K_F^{wgp}(X) \subseteq \mathcal O_F^{wgp}(X,\mathbb K(X))\) and \(\mathcal C_V^{wgp}(X) \subseteq \mathcal O_F^{wgp}(X,\mathbb F(X))\).
    In the case that \(\mathscr U\) is a countable family  of open sets, \(\mathscr U \in \mathcal O_F^{wgp}(X,\mathbb K(X)) \implies \mathscr U \in \mathcal K_F^{wgp}(X)\) and \(\mathscr U \in \mathcal O_F^{wgp}(X,\mathbb F(X)) \implies \mathscr U \in \mathcal C_V^{wgp}(X)\).
\end{proposition}
\begin{proof}
    The fact that \(\mathcal K_F^{wgp}(X) \subseteq \mathcal O_F^{wgp}(X,\mathbb K(X))\) follows immediately from Definition 5.5 of \cite{Li2016} and the observation that every \(\mathscr U \in \mathcal K_F^{wgp}(X)\) must be countable.

    To see that \(\mathcal C_V^{wgp}(X) \subseteq \mathcal O_F^{wgp}(X,\mathbb F(X))\), let \(\mathscr U \in \mathcal C_V^{wgp}(X)\) and let \(\varphi : \mathscr U \to \omega\) be the partition as in Definition 5.7 of \cite{Li2016}.
    Note that \(|\mathscr U| = \omega\).
    Let \(A \subseteq X\) be a proper closed sets and \(V_1, V_2 , \ldots, V_m\) be open so that \(V_j \setminus A \neq\emptyset\) for each \(j\).
    Set \(W_j = V_j \setminus A\) for each \(j\).
    Then, we can find \(n \in \omega\) and a finite set \(F \subseteq X\) so that \(\bigcap_{j=1}^m (X \setminus W_j) \subseteq \bigcup\varphi^{-1}(n)\), \(F \cap V_j \neq\emptyset\) for each \(j\), and \(F \cap \bigcup\varphi^{-1}(n) = \emptyset\).
    Since \(A \subseteq \bigcap_{j=1}^m (X \setminus W_j)\), we see that \(\mathscr U \in \mathcal O_F^{wgp}(X,\mathbb F(X))\).

    Let \(\mathscr U \in \mathcal O_F^{wgp}(X,\mathbb K(X))\) be countable and notice that Definition 5.5 of \cite{Li2016} is clearly satisfied.
    So let \(\mathscr U \in \mathcal O_F^{wgp}(X,\mathbb F(X))\) be countable.
    Let \(\varphi : \mathscr U \to \omega\) be as specified in the definition and \(V_1, V_2, \ldots, V_m\) be arbitrary proper non-empty open sets.
    Notice that \(A = \bigcap_{j=1}^m (X\setminus V_j)\) is a closed set.
    Since \((X\setminus A) \cap V_j \neq \emptyset\) for each \(j\), we can find \(n_0 \in \omega\) and a finite set \(F \subseteq X\) that satisfy the criteria listed above.
    This satisfies Definition 5.7 of \cite{Li2016}.
\end{proof}

The following proposition and examples show how the different cover types defined in this section relate to each other.

\begin{proposition}
    For a space \(X\) and an ideal \(\mathcal A\) consisting of closed sets,
    \[
        \begin{array}{ccccccccc}
            \Gamma(X,\mathcal A) & \subseteq & \mathcal O^{gp}(X, \mathcal A) & \subseteq & \mathcal O(X,\mathcal A) & \subseteq & \mathcal O^{wgp}(X, \mathcal A) & \subseteq & \mathcal O_X\\
            && \rotatebox{90}{\(\subseteq\)} && \rotatebox{90}{\(\subseteq\)} && \rotatebox{90}{\(\subseteq\)} &&\\
            && \mathcal O_F^{gp}(X, \mathcal A) & \subseteq & \mathcal O_F(X,\mathcal A) & \subseteq & \mathcal O_F^{wgp}(X, \mathcal A)
        \end{array}
    \]
\end{proposition}

\begin{example}
In general, \(\Gamma(X, \mathcal A)\) is a proper subset of \(\mathcal O^{gp}(X,\mathcal A)\).
\end{example}
Consider \(X = \mathbb R\) and \(\mathcal A = [\mathbb R]^{<\omega}\).
Let \(U_n = (-n-2,n+2)\) and \(V_n = (-2^{-n},2^{-n})\) for each \(n \in \omega\).
Notice that
\[
    \mathscr U = \{U_n : n \in \omega\} \cup \{V_n : n \in \omega\} \in \mathcal O^{gp}(X,\mathcal A) \setminus \Gamma(X, \mathcal A).
\]

\begin{example} \label{example:GroupProper}
In general, \(\mathcal O^{gp}(X, \mathcal A)\) is a proper subset of \(\mathcal O(X,\mathcal A)\).
\end{example}

Consider \(\omega_1\) with the discrete topology, \(\mathcal A\) the ideal of finite subsets, and \(\mathscr U\) consisting of precisely the set of finite subsets of \(\omega_1\).
Note that \(\mathscr U \in \mathcal O(\omega_1 , \mathcal A)\).
By way of contradiction, let \(\phi: \mathscr U \to \omega_1\) be so that \(\phi^{-1}(\alpha)\) is a finite collection of finite subsets for each \(\alpha < \omega_1\) that satisfies the groupability criterion.
Since \(\phi^{-1}(\alpha)\) is a finite family of finite sets, we assume that \(\phi^{-1}(\alpha)\) is a finite subset of \(\omega_1\).
For \(n\in\omega\), let \(\alpha_n < \omega_1\) be so that, for all \(\beta \geq \alpha_n\), \(\{ k : k < n \} \subseteq \phi^{-1}(\beta)\).
Now, let \(\alpha = \sup\{\alpha_n : n \in \omega\}\) and notice that \(\alpha < \omega_1\).
Since \(\alpha_n \leq \alpha\), we know that \(\{ k : k < n \} \subseteq \phi^{-1}(\alpha)\) for all \(n \in \omega\).
However, \(\phi^{-1}(\alpha)\) was assumed to be finite.
Therefore, \(\mathscr U \not\in \mathcal O^{gp}(\omega_1,\mathcal A)\).

\begin{example}
In general, \(\mathcal O(X, \mathcal A)\) is a proper subset of \(\mathcal O^{wgp}(X, \mathcal A)\).
\end{example}
Consider \(\mathbb R\) and \(\mathcal A = [\mathbb R]^{<\omega}\).
Let
\[
\mathcal U_n = \left\{ \left( \frac{j-1}{2^n}, \frac{j+1}{2^n} \right) : -3^n < j < 3^n \right\}.
\]
and set \(\mathcal U = \bigcup_n \mathcal U_n\).
Then \(\{-1,1\} \not\subseteq U\) for any \(U \in \mathcal U\), so \(\mathcal U\) is not an \(\omega\)-cover.
However, \(\left( -\left(\frac{3}{2}\right)^n,\left(\frac{3}{2}\right)^n\right) \subseteq \bigcup \mathcal U_n\) for each \(n\).
So \(\mathcal U\) is \(\omega\)-weakly groupable.

\begin{example}
In general, \(\mathcal O^{wgp}(X, \mathcal A)\) is a proper subset of \(\mathcal O_X\).
\end{example}
Consider \(\omega\) and let \(\mathscr U = \{\{n\} : n \in \omega\}\).
Let \(\varphi : \mathscr U \to \omega\) be arbitrary and consider \(\varphi^{-1}(0)\).
Then \(F = \bigcup \varphi^{-1}(0) \in [\omega]^{<\omega}\).
Let \(m = \max F\) and consider \(E = F \cup \{m+1\} \in [\omega]^{<\omega}\).
Clearly, \(E \not\subseteq \bigcup\varphi^{-1}(0)\).
Let \(n \in \omega\) be so that \(n > 0\) and notice that \(F \cap \bigcup\varphi^{-1}(n) = \emptyset\) so \(E \not\subseteq \bigcup\varphi^{-1}(n)\).
So \(\mathscr U \in \mathcal O_X \setminus \mathcal O^{wgp}(X, \mathcal A)\).

\subsection{Dense Sets and Blades}

\begin{notn}
    We will be using the following classes.
    \begin{itemize}
        \item
	    \(\Omega_{X,x}\) to be the set of all \(A \subseteq X\) with \(x \in \text{cl}_X(A)\).
	    We also call \(A \in \Omega_{X,x}\) a \emph{blade} of \(x\).
	    \item
	    \(\mathcal D_X\) to be the collection of all dense subsets of \(X\).
    \end{itemize}
\end{notn}

\begin{definition}
    For a space \(X\) and a dense \(D \subseteq X\), we say \(D\) is \emph{\(\omega\)-groupable} if there is a function \(\varphi:D \to |D|\) so that
    \begin{itemize}
        \item \(\varphi^{-1}(\xi)\) is finite for all \(\xi < |D|\) and
        \item for all open \(U \subseteq X\), there is a \(\xi_0 < |D|\) so that for all \(\xi \geq \xi_0\), \(U \cap \varphi^{-1}(\xi) \neq \emptyset\).
    \end{itemize}
    Let \(\mathcal D_X^{gp}\) denote the \(\omega\)-groupable dense subsets of \(X\).
\end{definition}
Observe that if \(D\) is countable, then the notion of \(\omega\)-groupable is equivalent to the standard notion of being groupable as discussed in \cite{KocinacEtAl2005}.

The following modification of denseness has shown up repeatedly in the literature without name: see \cite[Theorems 23 and 24]{KocinacEtAl2005} and \cite[Theorems 5.6 and 5.8]{Li2016}.
We give it a symbol here for notational convenience.
This definition does not seem to be a general topological notion, but rather is specific to the context of hyperspaces (or perhaps topologies on lattices).

\begin{definition}
    Fix a space \(X\).
    Consider the space of closed sets \(\mathbb F(X)\) with topology \(\mathscr T\).
    Then \(D \subseteq \mathbb F(X)\) is \emph{weakly \(\omega\)-groupable} if there is a function \(\varphi:D \to |D|\) so that
    \begin{itemize}
        \item \(\varphi^{-1}(\xi)\) is finite for all \(\xi < |D|\) and
        \item for all open \(V \subseteq \mathbb F(X)\), there is an \(\xi < |D|\) so that \(\bigcap \varphi^{-1}(\xi) \in V\).
    \end{itemize}
    Let \(\mathcal D_{\mathbb F(X)}^{wgp}\) denote the weakly \(\omega\)-groupable dense subsets of \(\mathbb F(X)\).
\end{definition}

\subsection{The Hyperspace Topologies}

The Vietoris and Fell hyperspace topologies are probably the most commonly studied hyperspace topologies.
Originally, the Fell topology, introduced in \cite{Fell}, was studied in the context of functional analysis and \(C^*\)-algebras.
For a detailed study of the Vietoris topology, see \cite{MichaelVietoris}.
We present a general version of these topologies in this section.

\begin{definition}
    Suppose \(X\) is a topological space and \(\mathcal A\) consists of closed sets.
    The \emph{upper \(\mathcal A\)-topology} or \emph{co-\(\mathcal A\)} topology on \(\mathbb F(X)\), denoted \(\mathscr T_{\mathcal A^+}\), is generated by subbasic open sets of the form
    \[
    (X \setminus A)^+ = \{F \in \mathbb F(X) : F \subseteq X \setminus A\} = \{F \in \mathbb F(X) : F \cap A = \emptyset\}
    \]
    along with \(\mathbb F(X)\). If \(\mathcal A\) is an ideal, then these are the basic open sets. Set \(\mathbb F(X,\mathcal A^+) = (\mathbb F(X), \mathscr T_{\mathcal A^+})\).
    \begin{itemize}
        \item
        The \emph{upper Fell topology} (or co-compact topology), \(\mathscr T_{F^+}\) is the case where \(\mathcal A\) is the collection of compact subsets of \(X\).
        Set \(\mathbb F(X,F^+) = (\mathbb F(X), \mathscr T_{F^+})\).
        \item
        The \emph{upper Vietoris topology}, \(\mathscr T_{V^+}\) is the case where \(\mathcal A\) is the collection of closed subsets of \(X\).
        Set \(\mathbb F(X,V^+) = (\mathbb F(X), \mathscr T_{V^+})\).
        \item
        The \emph{co-finite topology}, \(\mathscr T_{Z^+}\) is the case where \(\mathcal A\) is the collection of finite subsets of \(X\). Set \(\mathbb F(X,Z^+) = (\mathbb F(X), \mathscr T_{Z^+})\).
    \end{itemize}
\end{definition}

\begin{definition}
    Suppose \(X\) is a topological space.
    Let \(\mathcal A\) be an ideal consisting of closed subsets of \(X\).
    The topology on \(\mathbb F(X)\) generated by \(\mathcal A\), denoted \(\mathscr T_{\mathcal A}\), has basic open sets of the form
    \[
    [A; V_1,\ldots,V_n] = \{F \in \mathbb F(X) : F \cap A = \emptyset \text{ and } F \cap V_j \neq \emptyset \text{ for } 1 \leq j \leq n\}
    \]
    where \(A \in \mathcal A\) and \(V_1,\ldots,V_n \subseteq X\) are open. Set \(\mathbb F(X,\mathcal A) = (\mathbb F(X), \mathscr T_{\mathcal A})\).

    \begin{itemize}
    \item
        The \emph{Fell topology}, \(\mathscr T_F\) is the case where where \(\mathcal A\) is the collection of compact subsets of \(X\). Set \(\mathbb F(X,F) = (\mathbb F(X), \mathscr T_F)\).
        \item
        The \emph{Vietoris topology}, \(\mathscr T_V\) is the case where where \(\mathcal A\) is the collection of closed subsets of \(X\). Set \(\mathbb F(X,V) = (\mathbb F(X), \mathscr T_V)\).
    \end{itemize}
\end{definition}
One ought to notice that the Vietoris topology in this context is defined on the space of proper, non-empty closed sets and not just on the space of non-empty compact subsets of \(X\).

\begin{proposition} \label{prop:VietorisNonempty}
    Suppose \(X\) is a topological space.
    Let \(\mathcal A\) be an ideal consisting of closed subsets of \(X\).
    If the basic open set \([A; V_1, \ldots, V_n]\) in \(\mathbb F(X,\mathcal A)\) is non-empty, then \(V_j \cap (X \setminus A) \neq \emptyset\) for all \(j\).
\end{proposition}
\begin{proof}
    Let \(G \in [A; V_1, \ldots, V_n]\). Then \(G \subseteq X \setminus A\) and \(G \cap V_j \neq \emptyset\) for all \(j\).
    Thus \((X \setminus A) \cap V_j \neq \emptyset\) for all \(j\).
\end{proof}

\begin{lemma} \label{lem:IdealContainmentTopologyContainment}
    If \(\mathcal B \subseteq \mathcal A\) are ideals of closed sets in a space \(X\), then
    \begin{itemize}
        \item
        \(\mathscr T_{B^+} \subseteq \mathscr T_{A^+}\) and
        \item
        \(\mathscr T_{B} \subseteq \mathscr T_{A}\).
    \end{itemize}
\end{lemma}
\begin{proof}
    Consider an open set in \(\mathbb F(X, \mathcal B^+)\) which is necessarily of the form \((X \setminus B)^+\) for \(B \in \mathcal B\).
    Since \(B \in \mathcal A\), we see that \((X \setminus B)^+\) is open in \(\mathbb F(X, \mathcal A^+)\).

    Now suppose we have an open set \(U\) in \(\mathbb F(X, \mathcal B)\).
    Let \(F \in U\) be arbitrary and let \(B \in \mathcal B\) and \(V_1, \ldots, V_n\) open in \(X\) be so that
    \[
        F \in [B; V_1 , \ldots , V_n ] \subseteq U.
    \]
    Since \(B \in \mathcal A\), notice that \([B; V_1 , \ldots, V_n ]\) is a \(\mathscr T_{\mathcal A}\) neighborhood of \(F\) and, as \(F\) was chosen to be arbitrary, \(U\) is open in \(\mathbb F(X, \mathcal A)\).
\end{proof}

\subsection{Networks}

Recall the notion of a \(\pi\)-network.

\begin{definition}
    For a space \(X\) and \(\mathscr Y \subseteq \wp(X)\), \(\mathscr Y\) is said to be a \emph{\(\pi\)-network} if, for every non-empty open subset \(U\) of \(X\), there exists some \(Y \in \mathscr Y\) so that \(Y \subseteq U\).
\end{definition}

The following definition is motivated by the  \(\Pi_\omega\) and \(\Pi_k\) used in \cite{KocinacEtAl2005}.
\begin{definition}
    For a family \(\mathcal A \subseteq \wp(X)\), let \(\Pi_{\mathcal A}(X)\) denote the set of all \(\pi\)-networks consisting of sets from \(\mathcal A\).
\end{definition}

Note that
\begin{itemize}
    \item
    \(\Pi_k\) is \(\Pi_{\mathbb K(X)}\) and
    \item \(\Pi_\omega\) is \(\Pi_{[X]^{<\omega}}\).
\end{itemize}

In Section 3 of \cite{Li2016}, Li defined a version of \(\pi\)-network that is relevant to the full Vietoris and Fell topologies.
These are called \(\pi_V\)- and \(\pi_F\)-networks.
The definition below generalizes both of these.

\begin{definition} \label{def:piAF}
    Let \(\mathcal A\) be an ideal of closed subsets of \(X\) and let \(\mathscr X_{\mathcal A} \subseteq \mathcal A \times \mathscr T_X^{<\omega}\) be defined by the rule
    \[
        \langle A, V_1 , \ldots , V_n \rangle \in \mathscr X_{\mathcal A} \iff (\forall j)[V_j \setminus A \neq \emptyset].
    \]
    We say that a collection \(\mathscr Y \subseteq \mathscr X_{\mathcal A}\) is a \emph{\(\pi_{\mathcal A, F}\)-network} of \(X\) if, for each proper open subset \(U\) of \(X\), there exists \(\langle A , V_1 , \ldots , V_n \rangle \in \mathscr Y\) and a finite \(F \subseteq X\) so that \(F \cap V_j \neq \emptyset\) for each \(j\), \(A \subseteq U\), and \(F \cap U = \emptyset\).
    Let \(\Pi^F_{\mathcal A}(X)\) denote the family of all \(\pi_{\mathcal A,F}\)-networks of \(X\).
\end{definition}

\begin{lemma} \label{lem:IdealContainmentNetworkContainment}
    Suppose \(\mathcal B \subseteq \mathcal A\) are ideals of closed sets in \(X\).
    Then
    \begin{itemize}
        \item
        \(\Pi_{\mathcal B}(X) \subseteq \Pi_{\mathcal A}(X)\) and
        \item
        \(\Pi^F_{\mathcal B}(X) \subseteq \Pi^F_{\mathcal A}(X)\).
    \end{itemize}
\end{lemma}
\begin{proof}
    The first claim is clear.

    Now suppose \(\mathscr Y \in \Pi^F_{\mathcal B}(X)\) and let \(U\) be open in \(X\).
    Then there exists \(\langle B , V_1 , \ldots , V_n \rangle \in \mathscr Y\) and \(F \in [X]^{<\omega}\) that satisfy the conditions in the definition.
    Since \(B \in \mathcal A\), we see that \(\mathscr Y \in \Pi^F_{\mathcal A}(X)\).
\end{proof}

\begin{proposition} \label{prop:NetworksEquiv}
    Let \(X\) be a topological space and \(\mathcal A\) be an ideal of closed subsets of \(X\).
    \begin{itemize}
        \item
        \(\mathscr Y \in \Pi^F_{\mathbb K(X)}(X)\) if and only if \(\mathscr Y\) is a \(\pi_F\)-network,
        \item
        If \(\mathscr Y \in \Pi^F_{\mathbb F(X)}(X)\), then \(\{\langle X \setminus A, V_1, \ldots, V_n \rangle :  \langle A, V_1, \ldots, V_n \rangle \in \mathscr Y\}\) is a \(\pi_V\)-network, and
        \item If \(\mathscr Y\) is a a \(\pi_V\)-network, then \(\{\langle \bigcap_{j=1}^n (X \setminus V_j), V_1, \ldots, V_n \rangle : \langle V_1, \ldots, V_n \rangle \in \mathscr Y\} \in \Pi^F_{\mathbb F(X)}(X)\).
    \end{itemize}
\end{proposition}
\begin{proof}
    The first equivalence comes directly from the definitions in \cite{Li2016}.
    For the \(\pi_V\)-networks, first suppose that \(\mathscr Y \in \Pi^F_{\mathbb F(X)}(X)\).
    Let \(U \subseteq X\) be a proper open set.
    Then there is \(\langle A, V_1, \ldots, V_n \rangle \in \mathscr Y\) and \(F \subseteq X\) finite so that \(F \cap V_j \neq \emptyset\) (for each \(j\)), \(A \subseteq U\), and \(F \cap U = \emptyset\).
    Set \(V_0 = X \setminus A\). Then since \((X \setminus A) \subseteq \bigcup_{j=0}^n V_j\),
    \[
    \bigcap_{j=0}^n (X \setminus V_j) = X \setminus \left(\bigcup_{j=0}^n V_j \right) \subseteq A.
    \]
    So \(\{\langle X \setminus A, V_1, \ldots, V_n \rangle :  \langle A, V_1, \ldots, V_n \rangle \in \mathscr Y\}\) is a \(\pi_V\)-network.

    Now suppose \(\mathscr Y\) is a \(\pi_V\)-network.
    Let \(U \subseteq X\) be a proper open set.
    Then we can find \(\langle V_1,\ldots, V_n \rangle \in \mathscr Y\) and \(F \subseteq X\) finite so that
    \begin{itemize}
        \item \(F \cap V_j \neq \emptyset\) for \(1 \leq j \leq n\),
        \item \(\bigcap_{j=1}^n (X \setminus V_j) \subseteq U\), and
        \item \(F \cap U = \emptyset\).
    \end{itemize}
    This is enough to show that \(\{\langle \bigcap_{j=1}^n (X \setminus V_j), V_1, \ldots, V_n \rangle : \langle V_1, \ldots, V_n \rangle \in \mathscr Y\} \in \Pi^F_{\mathbb F(X)}(X)\).
\end{proof}

\subsection{Selection Principles and Games}

Though games of countable-length are more commonly studied, we extend our considerations to ordinal-length games for more generality.
Ordinal-games have been introduced before (see \cite{ScheepersTall}).

\begin{definition}
	Given a set \(\mathcal A\) and another set \(\mathcal B\), we define the \emph{finite selection game} \(G^\alpha_{\text{fin}}(\mathcal A, \mathcal B)\) for \(\mathcal A\) and \(\mathcal B\) as follows:
	\[
		\begin{array}{c|cccccc}
			\text{I} & A_0 & A_1 & A_2 & \ldots & A_\xi & \ldots\\
			\hline
			\text{II} & \mathcal F_0 & \mathcal F_1 & \mathcal F_2 & \ldots & \mathcal F_\xi & \ldots
		\end{array}
	\]
	where \(A_\xi \in \mathcal A\) and \(\mathcal F_\xi \in [A_\xi]^{<\omega}\) for all \(\xi < \alpha\).
	We declare Two the winner if \(\bigcup\{ \mathcal F_\xi : \xi < \alpha \} \in \mathcal B\).
	Otherwise, One wins.
	We let \(G_{\text{fin}}(\mathcal A, \mathcal B)\) denote \(G^\omega_{\text{fin}}(\mathcal A, \mathcal B)\).
\end{definition}

\begin{definition}
	Similarly, we define the \emph{single selection game} \(G^\alpha_1(\mathcal A, \mathcal B)\) as follows:
	\[
		\begin{array}{c|cccccc}
			\text{I} & A_0 & A_1 & A_2 & \ldots & A_\xi & \ldots\\
			\hline
			\text{II} & x_0 & x_1 & x_2 & \ldots & x_\xi & \ldots
		\end{array}
	\]
	where each \(A_\xi \in \mathcal A\) and \(x_\xi \in A_\xi\).
	We declare Two the winner if \(\{ x_\xi : \xi \in \alpha \} \in \mathcal B\).
	Otherwise, One wins.
	We let \(G_{1}(\mathcal A, \mathcal B)\) denote \(G^\omega_{1}(\mathcal A, \mathcal B)\).
\end{definition}

For collections \(\mathcal A\) and \(\mathcal B\) and an ordinal \(\alpha\), we refer to \(G^\alpha_1(\mathcal A, \mathcal B)\) and \(G^\alpha_{\text{fin}}(\mathcal A, \mathcal B)\) as \emph{selection games}.

\begin{definition}
    We define strategies of various strength below.
    \begin{itemize}
    \item A \emph{strategy for player One} in \(G^\alpha_1(\mathcal A, \mathcal B)\) is a function \(\sigma:(\bigcup \mathcal A)^{<\alpha} \to \mathcal A\).
    A strategy \(\sigma\) for One is called \emph{winning} if whenever \(x_\xi \in \sigma\langle x_\zeta : \zeta < \xi \rangle\) for all \(\xi < \alpha\), \(\{x_\xi:\xi \in \alpha\} \not\in \mathcal B\).
    If player One has a winning strategy, we write \(\One \uparrow G^\alpha_1(\mathcal A, \mathcal B)\).
    \item A \emph{strategy for player Two} in \(G^\alpha_1(\mathcal A, \mathcal B)\) is a function \(\tau:\mathcal A^{<\alpha} \to \bigcup \mathcal A\).
    A strategy \(\tau\) for Two is \emph{winning} if whenever \(A_\xi \in \mathcal A\) for all \(\xi < \alpha\), \(\{\tau(A_0,\ldots,A_\xi) : \xi < \alpha\} \in \mathcal B\).
    If player Two has a winning strategy, we write \(\Two \uparrow G^\alpha_1(\mathcal A, \mathcal B)\).
    \item A \emph{predetermined strategy} for One is a strategy which only considers the current turn number.
    We call this kind of strategy predetermined because One is not reacting to Two's moves, they are just running through a pre-planned script. Formally it is a function \(\sigma:\alpha \to \mathcal A\).
    If One has a winning predetermined strategy, we write \(\One \underset{\text{pre}}{\uparrow} G^\alpha_1(\mathcal A, \mathcal B)\).
    \item A \emph{Markov strategy} for Two is a strategy which only considers the most recent move of player One and the current turn number.
    Formally it is a function \(\tau:\mathcal A \times \alpha \to \bigcup \mathcal A\).
    If Two has a winning Markov strategy, we write \(\Two \underset{\text{mark}}{\uparrow} G^\alpha_1(\mathcal A, \mathcal B)\).
    \item
    If there is a single element \(x_0 \in \mathcal A\) so that the constant function with value \(x_0\) is a winning strategy for One, we say that One has a \emph{constant winning strategy}, denoted by \(\One \underset{\text{cnst}}{\uparrow} G^\alpha_1(\mathcal A, \mathcal B)\).
    \end{itemize}
    These definitions can be extended to \(G^\alpha_{\text{fin}}(\mathcal A, \mathcal B)\) in the obvious way.
\end{definition}

\begin{definition}
    Two games \(\mathcal G_1\) and \(\mathcal G_2\) are said to be \emph{strategically dual} provided that the following two statements hold:
    \begin{itemize}
        \item \(\text{I} \uparrow \mathcal G_1 \text{ iff } \text{II} \uparrow \mathcal G_2\)
        \item \(\text{I} \uparrow \mathcal G_2 \text{ iff } \text{II} \uparrow \mathcal G_1\)
    \end{itemize}
    Two games \(\mathcal G_1\) and \(\mathcal G_2\) are said to be \emph{Markov dual} provided that the following two statements hold:
    \begin{itemize}
        \item \(\text{I} \underset{\text{pre}}{\uparrow} \mathcal G_1 \text{ iff } \text{II} \underset{\text{mark}}{\uparrow} \mathcal G_2\)
        \item \(\text{I} \underset{\text{pre}}{\uparrow} \mathcal G_2 \text{ iff } \text{II} \underset{\text{mark}}{\uparrow} \mathcal G_1\)
    \end{itemize}
    Two games \(\mathcal G_1\) and \(\mathcal G_2\) are said to be \emph{dual} provided that they are both strategically dual and Markov dual.
\end{definition}

\begin{definition}
    The reader may be more familiar with selection principles than selection games.
    Let \(\mathcal A\) and \(\mathcal B\) be collections and \(\alpha\) be an ordinal.
    The selection principle \(S_1^\alpha(\mathcal A, \mathcal B)\) for a space \(X\) is the following property:
    Given any \(\alpha\)-length sequence \(\langle A_\beta : \beta < \alpha \rangle\) from \(\mathcal A\), there exists \(\{x_\beta : \beta < \alpha \}\) with \(x_\beta \in A_\beta\) for each \(\beta < \alpha\) so that \(\{ x_\beta : \beta < \alpha \} \in \mathcal B\).
    \(S_{\text{fin}}^\alpha(\mathcal A, \mathcal B)\) is similarly defined, but with finite selections instead of single selections.
\end{definition}

\begin{remark}
    In general, \(S_\square^\alpha(\mathcal A, \mathcal B)\) holds if and only if \(\text{I} \underset{\text{pre}}{\not\uparrow} G_\square^\alpha(\mathcal{A},\mathcal{B})\) where \(\square \in \{1 , \text{fin} \}\).
    See \cite[Prop. 15]{ClontzDuality}.
\end{remark}

\begin{definition}
    An even more fundamental type of selection is inspired by the Lindel\"{o}f property.
    Let \(\mathcal A\) and \(\mathcal B\) be collections.
    Then \({\mathcal A \choose \mathcal B}\) means that, for every \(A \in \mathcal A\), there exists \(B \subseteq A\) so that \(B \in \mathcal B\).
    Scheepers calls this a \emph{Bar-Ilan selection principle} in \cite{Scheepers2003}.
\end{definition}

\begin{remark} \label{rmk:ClontzBarIlan}
    Let \(\mathcal A\) and \(\mathcal B\) be collections, and \(\kappa\) be a cardinal.
    We let \(\mathcal B^\kappa = \{B \in \mathcal B: |B| \leq \kappa\}\).
    Then One fails to have a constant strategy in \(G_1^\kappa(\mathcal A, \mathcal B)\) if and only if \({\mathcal A \choose \mathcal B^\kappa}\) holds as shown in \cite[Prop. 15]{ClontzDuality}.
\end{remark}

\begin{remark}
    Note the following relationship between strategies in selection games and other types of selection principles.
    \begin{center}
        \begin{tikzpicture}
            \node [below] at (-6,1) {\(\Two\!\! \underset{\text{mark}}{\uparrow} \!\! G^\kappa_\square(\mathcal A, \mathcal B)\)};
            \node [below] at (-3,1) {\(\Two \uparrow G^\kappa_\square(\mathcal A, \mathcal B)\)};
            \node [below] at (0,1) {\(\One \not\uparrow G^\kappa_\square(\mathcal A, \mathcal B)\)};
            \node [below] at (3,1) {\(\One \underset{\text{pre}}{\not\uparrow} \mathcal G^\kappa_\square(\mathcal A, \mathcal B)\)};
            \node [below] at (6,1) {\(\One \underset{\text{cnst}}{\not\uparrow} G^\kappa_\square(\mathcal A, \mathcal B)\)};

            \foreach \x in {-4.5,-1.5,1.5,4.5} {
                \node at (\x,.65) {\(\Rightarrow\)};
            }
            \node [below] at (3,0) {\(\Updownarrow\)};
            \node [below] at (6,0) {\(\Updownarrow\)};
            \node at (4.5,-1.35) {\(\Rightarrow\)};
            \node [below] at (3,-1) {\(S_\square^\kappa(\mathcal A, \mathcal B)\)};
            \node [below] at (6,-1) {\({\mathcal A \choose \mathcal B^\kappa}\)};
        \end{tikzpicture}
    \end{center}
\end{remark}

\section{Translation Theorems} \label{section:Translations}

\begin{definition}
    For two selection games \(\mathcal G\) and \(\mathcal H\), we say that \(\mathcal G \leq_{\Two} \mathcal H\) if
    \begin{itemize}
        \item \(\Two \underset{\text{mark}}{\uparrow} \mathcal G \implies \Two \underset{\text{mark}}{\uparrow} \mathcal H\),
        \item \(\Two \uparrow \mathcal G \implies \Two \uparrow \mathcal H\),
        \item \(\One \not\uparrow \mathcal G \implies \One \not\uparrow \mathcal H\), and
        \item \(\One \not\underset{\text{pre}}{\uparrow} \mathcal G \implies \One \not\underset{\text{pre}}{\uparrow} \mathcal H\).
    \end{itemize}
    Note that \(\leq_{\Two}\) is transitive and that, if \(\mathcal G \leq_{\Two} \mathcal H\) and \(\mathcal H \leq_{\Two} \mathcal G\), then \(\mathcal G\) and \(\mathcal H\) are equivalent.
    In the case that \(\mathcal G \leq_{\Two} \mathcal H\) and \(\mathcal H \leq_{\Two} \mathcal G\), we write \(\mathcal G \equiv \mathcal H\).
\end{definition}

We start by recalling previous results about translating strategies.
\begin{theorem}[Theorem 12 of \cite{CCJH2020}]\label{OldTranslation}
    Let \(\mathcal A\), \(\mathcal B\), \(\mathcal C\), and \(\mathcal D\) be collections and \(\alpha\) be an ordinal.
    Suppose there are functions
    \begin{itemize}
        \item \(\overleftarrow{T}_{\One,\xi}:\mathcal B \to \mathcal A\) and
        \item \(\overrightarrow{T}_{\Two,\xi}: \left( \bigcup \mathcal A \right) \times \mathcal B \to \bigcup \mathcal B\)
    \end{itemize}
    for each \(\xi \in \alpha\), so that
    \begin{enumerate}[label=(Tr\arabic*)]
        \item \label{TranslationA} If \(x \in \overleftarrow{T}_{\One,\xi}(B)\), then \(\overrightarrow{T}_{\Two,\xi}(x,B) \in B\)
        \item \label{TranslationB} If \(x_\xi \in \overleftarrow{T}_{\One,\xi}(B_\xi)\) and \(\{x_\xi : \xi \in \alpha\} \in \mathcal C\), then \(\{\overrightarrow{T}_{\Two,\xi}(x_\xi,B_\xi) : \xi \in \alpha\} \in \mathcal D\).
    \end{enumerate}
    Then \(G^\alpha_1(\mathcal A,\mathcal C) \leq_{\Two} G^\alpha_1(\mathcal B, \mathcal D)\).
\end{theorem}

\begin{remark}
    Note that under the hypotheses of Theorem \ref{OldTranslation}, we also have that \({\mathcal A \choose \mathcal C^\alpha} \implies {\mathcal B \choose \mathcal D^\alpha}\).
    Set \(A = \overleftarrow{T}_{\One,0}(B)\), and then find \(C \subseteq A\) so that \(|C| \leq \alpha\) and \(C \in \mathcal C\), say \(C = \{x_\xi : \xi < \alpha\}\).
    Define \(D = \{\Ttwoxi(x_\xi,B) : \xi < \alpha\}\).
    Then \(|D| \leq \alpha\) and \(D \in \mathcal D\).
\end{remark}

\begin{corollary}[Corollary 13 of \cite{CCJH2020}]\label{corollary:EasyTranslate}
    Let \(\mathcal A\), \(\mathcal B\), \(\mathcal C\), and \(\mathcal D\) be collections and \(\alpha\) be an ordinal.
    Suppose there is a map \(\varphi : \left( \bigcup \mathcal B \right) \times \alpha \to \left( \bigcup \mathcal A \right)\) so that
    \begin{itemize}
        \item
        for all \(B \in \mathcal B\) and all \(\xi < \alpha\), \(\{ \varphi(y,\xi) : y \in B\} \in \mathcal A\)
        \item
        if \(\{ \varphi(y_\xi,\xi) : \xi < \alpha \} \in \mathcal C\), then \(\{ y_\xi : \xi < \alpha \} \in \mathcal D\)
    \end{itemize}
    Then \(G^\alpha_1(\mathcal A, \mathcal C) \leq_\Two G^\alpha_1(\mathcal B, \mathcal D)\).
\end{corollary}

We now provide the finite selection counterpart to Theorem \ref{OldTranslation}.
\begin{theorem}[The Translation Theorem] \label{TranslationFin}
    Let \(\mathcal A\), \(\mathcal B\), \(\mathcal C\), and \(\mathcal D\) be collections and let \(\alpha\) be an ordinal. Suppose there are functions
    \begin{itemize}
        \item \(\Tonexi:\mathcal B \to \mathcal A\) and
        \item \(\Ttwoxi:\left[\bigcup \mathcal A \right]^{<\omega} \times \mathcal B \to \left[\bigcup \mathcal B \right]^{<\omega}\)
    \end{itemize}
    for each \(\xi < \alpha\) so that
    \begin{enumerate}[label=(P\arabic*)]
        \item \label{FiniteTransA} If \(\mathcal F \in [\Tonexi(B)]^{<\omega}\), then \(\Ttwoxi(\mathcal F,B) \in [B]^{<\omega}\)
        \item \label{FiniteTransB} If \(\mathcal F_\xi \in [\Tonexi(B_\xi)]^{<\omega}\) for each \(\xi < \alpha\) and \(\bigcup_{\xi < \alpha} \mathcal F_\xi \in \mathcal C\), then \(\bigcup_{\xi < \alpha} \Ttwoxi(\mathcal F_\xi,B_\xi) \in \mathcal D\).
    \end{enumerate}
    Then \(G^\alpha_{\text{fin}}(\mathcal A, \mathcal C) \leq_{\Two} G^\alpha_{\text{fin}}(\mathcal B, \mathcal D)\).
\end{theorem}
\begin{proof}
    Suppose \(\tau\) is a winning Markov strategy for Two in \(G_{\text{fin}}^\alpha(\mathcal A, \mathcal C)\).
    Define a Markov strategy \(t\) for Two in \(G_{\text{fin}}^\alpha(\mathcal B, \mathcal D)\) by \(t(B,\xi) = \Ttwoxi(\tau(\Tonexi(B), \xi) , B)\).
    First note that if \(B \in \mathcal B\), then \(\tau(\Tonexi(B), \xi) \in  [\Tonexi(B)]^{<\omega}\), and so \(t(B,\xi) \in [B]^{<\omega}\).
    So \(t\) really is a Markov strategy.
    We now check that \(t\) is a winning strategy.
    Let \(\langle B_\xi : \xi < \alpha \rangle\) be a sequence from \(\mathcal B\).
    As \(\tau\) was a winning Markov strategy for Two in \(G_{\text{fin}}^\alpha(\mathcal A, \mathcal C)\), we have that \(\bigcup_{\xi < \alpha} (\tau(\Tonexi(B_\xi), \xi) \in \mathcal C\).
    Hence, \ref{FiniteTransB} asserts that \(\bigcup_{\xi < \alpha} t(B_\xi,\xi) \in \mathcal D\).
    So Two has a winning Markov strategy in \(G_{\text{fin}}^\alpha(\mathcal B, \mathcal D)\).

    Now Suppose \(\tau\) is a winning strategy for Two in \(G_{\text{fin}}^\alpha(\mathcal A, \mathcal C)\).
    We will define a winning strategy \(t\) for Two in \(G_{\text{fin}}^\alpha(\mathcal B, \mathcal D)\) recursively.

    Given \(B_0 \in \mathcal B\), let \(A_0 = \overset{\leftarrow}{T}_{\One,0}(B_0)\) and \(\mathcal F_0 = \tau(A_0) \in [A_0]^{<\omega}\).
    Then define \(t(B_0) = \mathcal G_0 = \overset{\rightarrow}{T}_{\Two,0}(\mathcal F_0,B_0)\).
    By \ref{FiniteTransA}, we see that \(t(B_0) \in [B_0]^{<\omega}\).

    For a given \(\beta < \alpha\), assume we have \(\langle A_\xi : \xi < \beta \rangle\) coming from \(\mathcal A\), \(\langle B_\xi : \xi < \beta \rangle\) coming from \(\mathcal B\), \(\langle \mathcal F_\xi : \xi < \beta \rangle\), and \(\langle \mathcal G_\xi : \xi < \beta \rangle\) all appropriately defined.
    Given \(B_\beta \in \mathcal B\), let \(A_\beta = \overset{\leftarrow}{T}_{\One, \beta}(B_\beta)\) and \(\mathcal F_\beta = \tau (A_0 , A_1 , \ldots , A_\beta) \in [A_\beta]^{<\omega}\).
    Define \(\mathcal G_\beta = \overset{\rightarrow}{T}_{\Two, \beta}(\mathcal F_\beta, B_\beta)\).
    Notice that \(\mathcal G_\beta \in [B_\beta]^{<\omega}\) by \ref{FiniteTransA}.
    So we define \(t(B_0, B_1, \ldots , B_\beta) = \mathcal G_\beta\).

    Now that \(t\) is defined, we must show that it is winning.
    Suppose we have a complete run of the game which involves \(\langle A_\xi : \xi < \alpha \rangle\), \(\langle B_\xi : \xi < \alpha \rangle\), \(\langle \mathcal F_\xi : \xi < \alpha \rangle\), and \(\langle \mathcal G_\xi : \xi < \alpha \rangle\) appropriately defined.
    Since \(\tau\) was assumed to be a winning strategy for Two in \(G_{\text{fin}}^\alpha(\mathcal A, \mathcal C)\), \(\bigcup_{\xi<\alpha}\mathcal F_\xi \in \mathcal C\).
    Hence, \ref{FiniteTransB} guarantees that \(\bigcup_{\xi < \alpha} \mathcal G_\xi \in \mathcal D\).
    That is, \(t\) is a winning strategy for Two in \(G_{\text{fin}}^\alpha(\mathcal B, \mathcal D)\).

    Suppose \(\sigma\) is a winning strategy for One in \(G^\alpha_{\text{fin}}(\mathcal B, \mathcal D)\).
    We define a strategy \(s\) for One in \(G^\alpha_{\text{fin}}(\mathcal A, \mathcal C)\) recursively as follows.
    Define \(B_0 = \sigma(\emptyset)\) and \(s(\emptyset) = A_0 = \Tonexi(\sigma(\emptyset))\).
    Then, for \(\mathcal F_0 \in [A_0]^{<\omega}\), set \(\mathcal G_0 = \overset{\rightarrow}{T}_{\Two,0}(\mathcal F_0, B_0)\).
    By \ref{FiniteTransA}, \(\mathcal G_0 \in [B_0]^{<\omega}\).

    Assume that, for \(\xi < \alpha\), we have \(\langle A_\eta:\eta < \xi\rangle\), \(\langle B_\eta:\eta < \xi\rangle\), \(\langle\mathcal F_\eta : \eta < \xi\rangle\), and \(\langle\mathcal G_\eta : \eta < \xi\rangle\) so that \(\mathcal F_\eta \in [A_\eta]^{<\omega}\) and \(\mathcal G_\eta = \overset{\rightarrow}{T}_{\Two,\eta}(\mathcal F_\eta, B_\eta)\) for all \(\eta < \xi\).
    By \ref{FiniteTransA}, \(\mathcal G_\eta \in [B_\eta]^{<\omega}\) for all \(\eta < \xi\).
    Then define \(B_{\xi} = \sigma(\langle \mathcal G_\eta : \eta < \xi \rangle)\), \(A_{\xi} = \Tonexi(B_{\xi})\) and \(s(\langle \mathcal F_\eta : \eta < \xi \rangle) = A_{\xi}\).
    Now we need only show that \(s\) is a winning strategy for One in \(G^\alpha_{\text{fin}}(\mathcal A, \mathcal C)\).
    Since \(\sigma\) is a winning strategy for One in \(G^\alpha_{\text{fin}}(\mathcal B, \mathcal D)\), we have that
    \[
        \bigcup_{\xi < \alpha} \mathcal G_\xi \not\in \mathcal D.
    \]
    Since \(\mathcal F_\xi \in [A_\xi]^{<\omega}\) and \(\mathcal G_\xi = \Ttwoxi(\mathcal F_\xi, B_\xi)\), it must be the case that
    \[
        \bigcup_{\xi < \alpha} \mathcal F_\xi \not\in \mathcal C
    \]
    by \ref{FiniteTransB}.

    Finally, suppose \(\sigma\) is a pre-determined winning strategy for One in \(G^\alpha_{\text{fin}}(\mathcal B, \mathcal D)\).
    Define \(s : \alpha \to \mathcal A\) by \(s(\xi) = \Tonexi(\sigma(\xi))\).
    To see that \(s\) is a winning strategy for One in \(G^\alpha_{\text{fin}}(\mathcal A, \mathcal C)\), consider \(\langle\mathcal F_\xi : \xi < \alpha\rangle\) where \(\mathcal F_\xi \in [s(\xi)]^{<\omega}\) for all \(\xi < \alpha\).
    Let \(\mathcal G_\xi = \Ttwoxi(\mathcal F_\xi, \sigma(\xi))\) and notice that \(\mathcal G_\xi \in [\sigma(\xi)]^{<\omega}\) by \ref{FiniteTransA}.
    Since \(\sigma\) is a winning strategy for One in \(G^\alpha_{\text{fin}}(\mathcal B, \mathcal D)\), we have that
    \[
        \bigcup_{\xi < \alpha} \mathcal G_\xi \not\in \mathcal D.
    \]
    Since \(\mathcal F_\xi \in [s(\xi)]^{<\omega}\) and \(\mathcal G_\xi = \Ttwoxi(\mathcal F_\xi, \sigma(\xi))\), it must be the case that
    \[
        \bigcup_{\xi < \alpha} \mathcal F_\xi \not\in \mathcal C
    \]
    by \ref{FiniteTransB}.
\end{proof}

\begin{corollary} \label{cor:SimplifiedTranslation}
    Let \(\mathcal A\), \(\mathcal B\), \(\mathcal C\), and \(\mathcal D\) be collections and \(\alpha\) be an ordinal.
    Suppose there are functions
    \begin{itemize}
        \item \(\Tonexi :\mathcal B \to \mathcal A\) and
        \item \(\Ttwoxi : \left(\bigcup \mathcal A \right) \times \mathcal B \to \bigcup \mathcal B\)
    \end{itemize}
    for each \(\xi < \alpha\) so that the following two properties hold.
    \begin{enumerate}[label={(Ft\arabic*)}]
        \item \label{translationPropI} If \(x \in \Tonexi(B)\), then \(\Ttwoxi(x,B) \in B\).
        \item \label{translationPropII} If \(\mathcal F_\xi \in \left[\Tonexi(B_\xi)\right]^{<\omega}\) and \(\bigcup_{\xi < \alpha} \mathcal F_\xi \in \mathcal C\), then \(\bigcup_{\xi < \alpha} \left\{ \Ttwoxi(x,B_\xi) : x \in \mathcal F_\xi \right\} \in \mathcal D\).
    \end{enumerate}
    Then \(G_{\text{fin}}^\alpha(\mathcal A, \mathcal C) \leq_{\Two} G_{\text{fin}}^\alpha(\mathcal B, \mathcal D)\) and \(G_1^\alpha(\mathcal A, \mathcal C) \leq_{\Two} G_1^\alpha(\mathcal B , \mathcal D)\).
\end{corollary}
\begin{proof}
    Suppose we have \(\Tonexi\) and \(\Ttwoxi\) as stated in the hypothesis.
    We define \(\Sonexi = \Tonexi\) and \(\Stwoxi:\left[\bigcup \mathcal A\right]^{<\omega} \times \mathcal B \to \left[\bigcup B \right]^{<\omega}\) by
    \[
    \Stwoxi(\{x_1,\ldots,x_n\},B) = \left\{ \Ttwoxi(x_1,B), \ldots, \Ttwoxi(x_n,B) \right\}.
    \]
    Suppose \(F \in \left[\Sonexi(B)\right]^{<\omega}\), say \(F = \{x_1, \ldots, x_n\}\).
    By \ref{translationPropI}, each \(\Ttwoxi(x_j,B) \in B\), so
    \[
    \Stwoxi(F,B) = \left\{ \Ttwoxi(x_1,B), \ldots, \Ttwoxi(x_n,B) \right\} \in [B]^{<\omega}.
    \]
    This establishes \ref{FiniteTransA}.

    Now suppose \(F_\xi \in \left[\Sonexi(B_\xi)\right]^{<\omega}\) and \(\bigcup_{\xi < \alpha} F_\xi \in \mathcal C\).
    Then \(F_\xi \in \left[\Tonexi(B_\xi)\right]^{<\omega}\), so
    \[
    \bigcup_{\xi < \alpha} \Stwoxi(F_\xi,B_\xi) = \bigcup_{\xi < \alpha} \left\{ \Ttwoxi(x,B_\xi) : x \in \mathcal F_\xi \right\} \in \mathcal D
    \]
    by \ref{translationPropII}.
    This establishes \ref{FiniteTransB}.
    So \(\Sonexi\) and \(\Stwoxi\) witness that \(G_{\text{fin}}^\alpha(\mathcal A, \mathcal C) \leq_{\Two} G_{\text{fin}}^\alpha(\mathcal B, \mathcal D)\) by Theorem \ref{TranslationFin}.

    Now we need to establish that \(G_1^\alpha(\mathcal A, \mathcal C) \leq_{\Two} G_1^\alpha(\mathcal B , \mathcal D)\).
    First, notice that \ref{translationPropI} is identical to \ref{TranslationA}.
    To see that \ref{TranslationB} is satisfied, let \(x_\xi \in \Tonexi(B_\xi)\) for each \(\xi < \alpha\) be so that \(\{ x_\xi : \xi < \alpha \} \in \mathcal C\).
    Let \(\mathcal F_\xi = \{ x_\xi \}\) and notice that \(\mathcal F_\xi \in \left[ \bigcup \Tonexi(B_\xi) \right]^{<\omega}\) for each \(\xi < \alpha\).
    Moreover,
    \[
        \bigcup_{\xi < \alpha} \mathcal F_\xi = \{ x_\xi : \xi < \alpha \} \in \mathcal C.
    \]
    So, by \ref{translationPropII}, we see that
    \[
        \bigcup_{\xi < \alpha} \left\{ \Ttwoxi(x,B_\xi) : x \in \mathcal F_\xi \right\} = \left\{ \Ttwoxi(x_\xi,B_\xi) : \xi < \alpha \right\} \in \mathcal D.
    \]
    That is, we can apply Theorem \ref{OldTranslation} to conclude that \(G_1^\alpha(\mathcal A, \mathcal C) \leq_{\Two} G_1^\alpha(\mathcal B , \mathcal D)\).
\end{proof}

The following result is similar to Corollary \ref{corollary:EasyTranslate}.
\begin{corollary} \label{cor:Translation}
    Let \(\mathcal A\), \(\mathcal B\), \(\mathcal C\), and \(\mathcal D\) be collections and \(\alpha\) be an ordinal.
    Suppose there is a map \(\varphi:\left[\bigcup \mathcal B\right] \times \alpha \to \bigcup \mathcal A\) so that the following two conditions hold.
    \begin{itemize}
        \item For all \(B \in \mathcal B\) and all \(\xi \in \alpha\), \(\{\varphi(y,\xi) : y \in B\} \in \mathcal A\).
        \item If \(\mathcal G_\xi \in [B_\xi]^{<\omega}\) where \(B_\xi \in \mathcal B\) for each \(\xi < \alpha\) and \(\bigcup_{\xi < \alpha} \varphi[\mathcal G_\xi \times \{\xi\}] \in \mathcal C\), then \(\bigcup_{\xi < \alpha} \mathcal G_\xi \in \mathcal D\).
    \end{itemize}
    Then \(G^\alpha_{\text{fin}}(\mathcal A, \mathcal C) \leq_{\Two} G^\alpha_{\text{fin}}(\mathcal B, \mathcal D)\) and \(G_1^\alpha(\mathcal A, \mathcal C) \leq_{\Two} G_1^\alpha(\mathcal B , \mathcal D)\).
\end{corollary}
\begin{proof}
    Define \(\Tonexi :\mathcal B \to \mathcal A\) by
    \[
        \Tonexi(B) = \varphi[B \times \{\xi\}] = \{ \varphi(y,\xi) : y \in B\}.
    \]
    Then \(\Tonexi(B) \in \mathcal A\).

    We now define \(\Ttwoxi : \left[\bigcup \mathcal A \right] \times \mathcal B \to \bigcup \mathcal B\).
    Fix \(B\in\mathcal B\) and, for \(x \in \varphi[B \times \{\xi\}]\), choose \(y_{x,\xi} \in B\) so that \(\varphi(y_{x,\xi},\xi) = x\).
    Then set \(\Ttwoxi(x,B) = y_{x , \xi}\).
    For \(x \not\in \varphi[B \times \{\xi\}]\), set \(\Ttwoxi(x,B)\) to be an arbitrary element of \(\bigcup \mathcal B\).
    This guarantees that, if \(x \in \Tonexi(B)\), then \(\Ttwoxi(x,B) \in B\) which establishes \ref{translationPropI}.

    To check that \ref{translationPropII} holds, suppose \(\mathcal F_\xi \in \left[\Tonexi(B_\xi)\right]^{<\omega}\) and \(\bigcup_{\xi < \alpha} \mathcal F_\xi \in \mathcal C\).
    Let
    \[
        \mathcal G_\xi = \left\{ \Ttwoxi(x,B_\xi) : x \in \mathcal F_\xi \right\}
    \]
    Notice that
    \[
        x \in \mathcal F_\xi \implies x \in \Tonexi(B_\xi) \implies \mathcal G_\xi \in [B_\xi]^{<\omega}.
    \]
    Moreover, \(\varphi[\mathcal G_\xi\times\{\xi\}] = \mathcal F_\xi\) since \(\varphi\left(\Ttwoxi(x,B_\xi) , \xi\right) = x\) for each \(x \in \mathcal F_\xi\).
    That is, we have that
    \[
        \bigcup_{\xi < \alpha} \varphi[\mathcal G_\xi\times\{\xi\}] = \bigcup_{\xi < \alpha} \mathcal F_\xi \in \mathcal C
    \]
    which, by the hypotheses, allows us to conclude that
    \[
        \bigcup_{\xi < \alpha} \left\{ \Ttwoxi(x,B_\xi) : x \in \mathcal F_\xi \right\} = \bigcup_{\xi < \alpha} \mathcal G_\xi \in \mathcal D.
    \]
    Hence, Corollary \ref{cor:SimplifiedTranslation} applies and the proof is finished.
\end{proof}

Using these general translation theorems, we infer an easy to state, easy to apply, but less general version of the translation theorem.
This version is the main tool we reference later to prove results about selection games on the Vietoris and Fell hyperspace topologies.

\begin{corollary} \label{cor:CorrespondenceEquivalence}
    Suppose that \(\mathcal{A, B, C, D}\) are collections so that \(\bigcup \mathcal C \subseteq \bigcup \mathcal A\) and \(\bigcup \mathcal D \subseteq \bigcup \mathcal B\).
    Additionally, suppose that there exists a bijection \(\beta : \bigcup \mathcal A \to \bigcup \mathcal B\) with the following features:
    \begin{itemize}
        \item \(A \in \mathcal A \iff \beta[A] \in \mathcal B\), and
        \item \(C \in \mathcal C \iff \beta[C] \in \mathcal D\).
    \end{itemize}
    Let \(\alpha\) be an ordinal.
    Then \(G_1^\alpha(\mathcal A, \mathcal C) \equiv G^\alpha_1(\mathcal B, \mathcal D)\) and \(G^\alpha_{\text{fin}}(\mathcal A, \mathcal C) \equiv G^\alpha_{\text{fin}}(\mathcal B, \mathcal D)\).
\end{corollary}
\begin{proof}
    Note that, as \(A \in \mathcal A \iff \beta[A] \in \mathcal B\), we also have \(\beta^{-1}[B] \in \mathcal A \iff B \in \mathcal B\).
    Similarly, \(\beta^{-1}[D] \in \mathcal C \iff D \in \mathcal D\).
    By symmetry, we need only show that \(G^\alpha_1(\mathcal A, \mathcal C) \leq_\Two G^\alpha_1(\mathcal B, \mathcal D)\) and that \(G^\alpha_{\text{fin}}(\mathcal A, \mathcal C) \leq_\Two G^\alpha_{\text{fin}}(\mathcal B, \mathcal D)\).

    Define \(\varphi : \left(\bigcup \mathcal B\right) \times \alpha \to \bigcup \mathcal A\) by the rule \(\varphi(y,\xi) = \beta^{-1}(y)\).
    This is well-defined since \(\beta\) is a bijection.
    Let \(B \in \mathcal B\) and \(\xi < \alpha\) be fixed.
    Then notice that \[\{\varphi(y,\xi) : y \in B\} = \{ \beta^{-1}(y) : y \in B\} = \beta^{-1}[B] \in \mathcal A.\]
    Suppose \(\{\varphi(y_\xi,\xi) : \xi < \alpha\} \in \mathcal C\).
    That is, \(\{\beta^{-1}(y_\xi) : \xi < \alpha \} \in \mathcal C\) and so \(\{ y_\xi : \xi < \alpha \} \in \mathcal D\).
    Thus, by Corollary \ref{corollary:EasyTranslate}, \(G^\alpha_1(\mathcal A, \mathcal C) \leq_\Two G^\alpha_1(\mathcal B, \mathcal D)\).

    Suppose \(\mathcal G_\xi \in [B_\xi]^{<\omega}\) and \(\bigcup_{\xi<\alpha} \varphi[\mathcal G_\xi \times \{\xi\}] \in \mathcal C\).
    Notice that
    \[
        \bigcup_{\xi < \alpha} \varphi[\mathcal G_\xi \times \{\xi\}] = \bigcup_{\xi < \alpha} \beta^{-1}[\mathcal G_\xi] = \beta^{-1}\left[ \bigcup_{\xi < \alpha} \mathcal G_\xi \right] \in \mathcal C,
    \]
    so \(\bigcup_{\xi < \alpha} \mathcal G_\xi \in \mathcal D\).
    Thus, by Corollary \ref{cor:Translation}, \(G^\alpha_{\text{fin}}(\mathcal A, \mathcal C) \leq_\Two G^\alpha_{\text{fin}}(\mathcal B, \mathcal D)\).
\end{proof}

\begin{remark} \label{rmk:BarIlanThings}
    Under the assumptions of Corollary \ref{cor:CorrespondenceEquivalence}, one can prove that \({\mathcal A \choose \mathcal C^\alpha} \iff {\mathcal B \choose \mathcal D^\alpha}\) for all \(\alpha\) in a similar way to what was done above.
\end{remark}

\section{Equivalences with Covers} \label{section:Covers}

\subsection{The Upper Topologies}

\begin{lemma} \label{lemma:ACoverUpperFellDense}
    Let \(\mathcal A\) be an ideal of closed subsets of \(X\).
    \(\mathscr U \in \mathcal O(X,\mathcal A)\) if and only if \(c.\mathscr U \in \mathcal D_{\mathbb F(X,\mathcal A^+)}\).
\end{lemma}
\begin{proof}
    Let \(\mathscr U \in \mathcal O(X,\mathcal A)\) and let \((X \setminus A)^+\) be an arbitrary open set in \(\mathbb F(X,\mathcal A^+)\).
    Then, since \(\mathscr U\) is an \(\mathcal A\)-cover, there exists \(U \in \mathscr U\) so that \(A \subseteq U\).
    This is equivalent to \(X \setminus U \subseteq X \setminus A\), which establishes that \(c.\mathscr U \in \mathcal D_{\mathbb F(X,\mathcal A^+)}\).

    Now, suppose \(\{ X \setminus U : U \in \mathscr U \}\) is dense in the co-\(\mathcal A\) hyperspace.
    As in the situation above, consider any \(X \setminus A\) where \(A \in \mathcal A\).
    There exists some \(U\in\mathscr U\) so that \(X\setminus U \subseteq X \setminus A\) by density, which means \(A \subseteq U\).
    Hence, \(\mathscr U\) is an \(\mathcal A\)-cover of \(X\).
\end{proof}

\begin{lemma} \label{lemma:ACoverUpperFellBlade}
    Let \(\mathcal A\) be an ideal of closed subsets of \(X\) and \(G \in \mathbb F(X)\).
    Then \(\mathscr U \in \mathcal O(X, X \setminus G, \mathcal A)\) if and only if \(c.\mathscr U \in \Omega_{\mathbb F(X,\mathcal A^+),G}\).
\end{lemma}
\begin{proof}
    First suppose \(\mathscr U\) is an \(\mathcal A\)-cover of \(X \setminus G\) by open sets from \(X\).
    Let \((X \setminus A)^+\) be an arbitrary open neighborhood of \(G\) in \(\mathbb F(X,\mathcal A^+)\).
    Then \(G \subseteq X \setminus A\) and so \(A \subseteq X \setminus G\).
    Thus there exists \(U \in \mathscr U\) so that \(A \subseteq U\).
    So \((X \setminus U) \subseteq (X \setminus A)\).
    Now \(X \setminus U \in c.\mathscr U \cap (X \setminus A)^+\).
    Therefore \(c.\mathscr U \in \Omega_{\mathbb F(X,\mathcal A^+),G}\).

    Now, suppose \(c.\mathscr U \in \Omega_{\mathbb F(X,\mathcal A^+),G}\).
    Let \(A \in \mathcal A\) be so that \(A \subseteq X \setminus G\).
    Then \(G \subseteq X \setminus A\) and so \((X \setminus A)^+\) is an open neighborhood of \(G\).
    So there exists some \(U\in\mathscr U\) so that \(X \setminus U \subseteq X \setminus A\), which means \(A \subseteq U\).
    Hence, \(\mathscr U\) is a \(\mathcal A\)-cover of \(X \setminus G\) by open sets from \(X\).
\end{proof}

\begin{lemma} \label{lemma:ACoverUpperFellGroup}
    Let \(\mathcal A\) be an ideal of closed subsets of \(X\).
    \(\mathscr U \in \mathcal O^{gp}(X, \mathcal A)\) if and only if \(c.\mathscr U \in \mathcal D^{gp}_{\mathbb F(X,\mathcal A^+)}\).
\end{lemma}
\begin{proof}
    First suppose \(\mathscr U\) is an \(\omega\)-groupable \(\mathcal A\)-cover of \(X\) and set \(\mathscr D = c.\mathscr U\).
    Let \(\varphi:\mathscr U \to |\mathscr U|\) be a partition of \(\mathscr U\) as in the definition of \(\omega\)-groupable.
    Note that \(|\mathscr D| = |\mathscr U|\).
    Define \(\psi:\mathscr D \to |\mathscr U|\) by \(\psi(F) = \varphi(X \setminus F)\).
    Notice that \(X \setminus U \in \psi^{-1}(\xi)\) if and only if \(U \in \varphi^{-1}(\xi)\).
    Let \((X \setminus A)^+\) be an arbitrary open set in \(\mathbb F(X,\mathcal A^+)\).
    Thus there exists a \(\xi_0 < |\mathscr U|\) so that for all \(\xi \geq \xi_0\), there is a \(U_\xi \in \varphi^{-1}(\xi)\) so that \(A \subseteq U_\xi\).
    So \(X \setminus U_\xi \in \psi^{-1}(\xi)\) and \((X \setminus U_\xi) \subseteq (X \setminus A)\).
    Therefore \(c.\mathscr U \in \mathcal D^{gp}_{\mathbb F(X,\mathcal A^+)}\).

    Now, suppose \(c.\mathscr U \in \mathcal D^{gp}_{\mathbb F(X,\mathcal A^+)}\) and set \(\mathscr D = c.\mathscr U\).
    Let \(\varphi:\mathscr D \to |\mathscr D|\) be a partition of \(\mathscr D\) as in the definition of \(\omega\)-groupable.
    Note that \(|\mathscr U| = |\mathscr D|\).
    Define \(\psi:\mathscr U \to |\mathscr D|\) by \(\psi(U) = \varphi(X \setminus U)\).
    Notice that \(F \in \psi^{-1}(\xi)\) if and only if \(X \setminus F \in \varphi^{-1}(\xi)\).
    Let \(A \in \mathcal A\) and consider \((X \setminus A)^+\).
    There is a \(\xi_0 < |\mathscr D|\) so that for all \(\xi \geq \xi_0\) there is an \(F_\xi \in \varphi^{-1}(\xi)\) so that \(F_\xi \subseteq X \setminus A\).
    Set \(U_\xi = X \setminus F_\xi\).
    Then \(U_\xi \in \psi^{-1}(\xi)\) and \(A \subseteq U_\xi\).
    Hence, \(\mathscr U\) is an \(\omega\)-groupable \(\mathcal A\)-cover of \(X\).
\end{proof}

\begin{lemma} \label{lemma:ACoverUpperFellWeakGroup}
    Let \(\mathcal A\) be an ideal of closed subsets of \(X\).
    Then \(\mathscr U \in \mathcal O^{wgp}(X, \mathcal A)\) if and only if \(c.\mathscr U \in \mathcal D^{wgp}_{\mathbb F(X,\mathcal A^+)}\).
\end{lemma}
\begin{proof}
    First suppose \(\mathscr U\) is a weakly \(\omega\)-groupable \(\mathcal A\)-cover of \(X\) and set \(\mathscr D = c.\mathscr U\).
    Let \(\varphi:\mathscr U \to |\mathscr U|\) be a partition of \(\mathscr U\) as in the definition of weakly \(\omega\)-groupable.
    Note that \(|\mathscr D| = |\mathscr U|\).
    Define \(\psi:\mathscr D \to |\mathscr U|\) by \(\psi(F) = \varphi(X \setminus F)\).
    Notice that \(X \setminus U \in \psi^{-1}(\xi)\) if and only if \(U \in \varphi^{-1}(\xi)\).
    Let \(V = (X \setminus A)^+\) be an arbitrary open set in \(\mathbb F(X,\mathcal A^+)\).
    Thus there exists a \(\xi < |\mathscr U|\) so that \(A \subseteq \bigcup \varphi^{-1}(\xi)\).
    So \(\psi^{-1}(\xi) = \{X \setminus U : U \in \varphi^{-1}(\xi)\}\) and
    \[
   \bigcap \psi^{-1}(\xi) =  X \setminus \bigcup \varphi^{-1}(\xi) \subseteq X \setminus A
    \]
    Hence, \(\bigcap \psi^{-1}(\xi) \in V\).
    Therefore \(c.\mathscr U \in \mathcal D^{wgp}_{\mathbb F(X,\mathcal A^+)}\).

    Now, suppose \(c.\mathscr U \in \mathcal D^{wgp}_{\mathbb F(X,\mathcal A^+)}\) and set \(\mathscr D = c.\mathscr U\).
    Let \(\varphi:\mathscr D \to |\mathscr D|\) be a partition of \(\mathscr D\) as in the definition of weakly \(\omega\)-groupable.
    Note that \(|\mathscr U| = |\mathscr D|\).
    Define \(\psi:\mathscr U \to |\mathscr D|\) by \(\psi(U) = \varphi(X \setminus U)\).
    Notice that \(F \in \psi^{-1}(\xi)\) if and only if \(X \setminus F \in \varphi^{-1}(\xi)\).
    Let \(A \in \mathcal A\).
    Then \(V = (X \setminus A)^+\) is an open set in \(\mathbb F(X,\mathcal A^+)\).
    So there is a \(\xi < |\mathscr D|\) so that so that \(\bigcap \psi^{-1}(\xi) \in V\).
    Therefore
    \[
    \bigcap \psi^{-1}(\xi) \subseteq X \setminus A \implies A \subseteq \bigcup \varphi^{-1}(\xi).
    \]
    Hence, \(\mathscr U\) is a weakly \(\omega\)-groupable \(\mathcal A\)-cover of \(X\).
\end{proof}

\begin{lemma} \label{lem:UpperTopPiNetwork}
    Fix a space \(X\) and a collection \(\mathscr S = \{ (X \setminus A)^+ : A \in S\}\) of open subsets of \(\mathbb F(X,\mathcal A^+)\).
    Then \(\mathscr S \in \mathcal O_{\mathbb F(X,\mathcal A^+)}\) if and only if \(S \in \Pi_{\mathcal A}(X)\).
\end{lemma}
\begin{proof}
    First suppose \(\mathscr S \in \mathcal O_{\mathbb F(X,\mathcal A^+)}\).
    Let \(U \subseteq X\) be open.
    Then there is an \(A \in S\) so that \(X \setminus U \in (X \setminus A)^+\).
    Thus \(X \setminus U \subseteq X \setminus A\) and so \(A \subseteq U\).
    Therefore, \(S \in \Pi_{\mathcal A}(X)\).

    Now suppose that \(S \in \Pi_{\mathcal A}(X)\) and let \(F \in \mathbb F(X,\mathcal A^+)\).
    Then \(X \setminus F\) is open, so there is an \(A \in S\) so that \(A \subseteq X \setminus F\).
    Thus \(F \subseteq X \setminus A\) and so \(F \in (X \setminus A)^+\).
    Hence, \(\mathscr S \in \mathcal O_{\mathbb F(X,\mathcal A^+)}\).
\end{proof}

\subsection{The Full Topologies}

\begin{lemma} \label{lem:AFDensity}
    For any space \(X\) and any ideal \(\mathcal A\) consisting of closed sets, \(\mathscr U \in \mathcal O_F(X,\mathcal A)\) if and only if \(c.\mathscr U \in \mathcal D_{\mathbb F(X,\mathcal A)}\).
\end{lemma}
\begin{proof}
    Suppose \(\mathscr U\) is a \(\mathcal A_F\)-cover of \(X\).
    Let \([A;V_1,\ldots,V_n]\) be a non-empty basic open set in \(\mathbb F(X,\mathcal A)\).
    Then by Proposition \ref{prop:VietorisNonempty}, \(V_j \cap (X \setminus A) \neq \emptyset\) for all \(j\).
    Since \(\mathscr U\) is a \(\mathcal A_F\)-cover of \(X\), we can find \(U \in \mathscr U\) and \(E \subseteq X\) finite so that \(A \subseteq U\), \(E \cap V_j \neq\emptyset\) for each \(j\), and \(E \cap U = \emptyset\).
    Set \(F = X \setminus U\) and notice that \(F \in c.\mathscr U\), \(F \subseteq X \setminus A\), and since \(E \cap (X \setminus A) = \emptyset\), \(E \subseteq F\).
    Hence, \(\emptyset \neq E \cap V_j \subseteq F \cap V_j\) for each \(j\).
    That is, \(F \in [ A ; V_0 , \ldots , V_n]\), so \(c.\mathscr U\) is dense.

    Now suppose \(\mathscr D = c.\mathscr U\) is dense in \(\mathbb F(X,\mathcal A)\).
    Let \(A \in \mathcal A\) and \(V_1 , V_2 , \ldots, V_n \subseteq X\) be open so that \(V_j \cap (X \setminus A) \neq\emptyset\) for each \(j\).
    Notice that the basic open set \([A;V_1, \ldots, V_n]\) is non-empty, so since \(\mathscr D\) is dense, we can find \(F \in \mathscr D \cap [A;V_1, \ldots, V_n]\).
    Then \(F \subseteq X \setminus A\) and \(F \cap V_j \neq\emptyset\) for each \(j\).
    Immediately, we see that \(A \subseteq X \setminus F\).
    Set \(U = X \setminus F\) and notice that \(U \in \mathscr U\).
    Let \(x_j \in F \cap V_j\) and \(E = \{ x_1 , x_2 , \ldots , x_n \}\).
    Clearly, \(E \cap V_j \neq \emptyset\) for each \(j\) and, since \(E \subseteq F\), \(E \cap U = \emptyset\).
    That is, \(\mathscr U\) is a \(\mathcal A_F\)-cover of \(X\).
\end{proof}

\begin{lemma} \label{lem:AFcoversBlades}
    For any space \(X\), any closed set \(G \in \mathbb F(X)\), and any ideal \(\mathcal A\) consisting of closed sets, \(\mathscr U \in \mathcal O_F(X , X \setminus G, \mathcal A)\) if and only if \(c.\mathscr U \in \Omega_{\mathbb F(X, \mathcal A), G}\).
\end{lemma}
\begin{proof}
    Suppose \(\mathscr U \in \mathcal O_F(X , X \setminus G,\mathcal A)\) and let \([A;V_1, \ldots, V_n]\) be a basic neighborhood of \(G\) in \(\mathbb F(X,\mathcal A)\).
    Since \(G \subseteq X \setminus A\), \(A \subseteq X \setminus G\).
    Also \(V_j \cap G \neq \emptyset\) for each \(j\).
    Since \(\mathscr U\) is a \(\mathcal A_F\)-cover of \(X \setminus G\), we can find \(U \in \mathscr U\) and \(E \subseteq X\) finite so that \(A \subseteq U\), \(E \cap V_j \neq\emptyset\) for each \(j\), and \(E \cap U = \emptyset\).
    Set \(F = X \setminus U\) and notice that \(F \in c.\mathscr U\).
    Immediately, \(F \subseteq X \setminus A\).
    Moreover, since \(E \cap U = \emptyset\), we see that \(E \subseteq F\).
    Hence, for each \(j\), \(\emptyset \neq E \cap V_j \subseteq F \cap V_j\).
    Therefore, \(F \in [A;V_1,\ldots,V_n]\), and so \(G\) is in the closure of \(c.\mathscr U\).

    Now suppose \(c.\mathscr U \in \Omega_{\mathbb F(X, \mathcal A), G}\).
    Consider \(A \in \mathcal A\) with \(A \subseteq (X \setminus G)\) and \(V_1, \ldots, V_n \subseteq X\) open with \(V_j \cap (X \setminus (X \setminus G)) = V_j \cap G \neq \emptyset\) for all \(j\).
    Notice that \(G \in [A;V_1, \ldots, V_n]\).
    Thus we can find \(F \in c.\mathscr U\) so that both \(F \subseteq X \setminus A\) and \(F \cap V_j \neq\emptyset\) for each \(j\).
    Immediately, we see that \(A \subseteq X \setminus F\).
    Set \(U = X \setminus F\) and notice that \(U \in \mathscr U\).
    For each \(j\), choose \(x_j \in F \cap V_j\).
    Then let \(E = \{ x_1, \ldots, x_n\}\).
    Notice that \(E\) is a finite set with \(E \cap V_j \neq \emptyset\) for each \(j\).
    Moreover, \(E \subseteq F\) which means that \(E \cap U = \emptyset\).
    Therefore \(\mathscr U \in \mathcal O_F(X , X \setminus G,\mathcal A)\).
\end{proof}

\begin{lemma} \label{lem:AFGrp}
    For any space \(X\) and any ideal \(\mathcal A\) consisting of closed sets, \(\mathscr U \in \mathcal O_F^{gp}(X,\mathcal A)\) if and only if \(c.\mathscr U \in \mathcal D_{\mathbb F(X,\mathcal A)}^{gp}\).
\end{lemma}
\begin{proof}
    Suppose \(\mathscr U \in \mathcal O_F^{gp}(X,\mathcal A)\).
    Let \(\varphi:\mathscr U \to |\mathscr U|\) be the partition of \(\mathscr U\) into finite sets as specified in the definition of \(\mathcal O_F^{gp}(X,\mathcal A)\).
    Set \(\mathscr D = c.\mathscr U\), notice that \(|\mathscr D| = |\mathscr U|\), and define \(\psi:\mathscr D \to |\mathscr D|\) by \(\psi(F) = \varphi(X \setminus F)\).
    Notice that this partitions \(\mathscr D\) into finite sets and \(F \in \psi^{-1}(\xi)\) if and only if \(X \setminus F \in \varphi^{-1}(\xi)\).
    Let \([A;V_1,\ldots,V_m]\) be a non-empty basic open set in \(\mathbb F(X,\mathcal A)\). Then by Proposition \ref{prop:VietorisNonempty}, \(V_j \cap (X \setminus A) \neq \emptyset\) for all \(j\).
    Thus we can find a \(\xi_0 < |\mathscr D|\) so that for all \(\xi \geq \xi_0\), there are \(U_\xi \in \varphi^{-1}(\xi)\) and \(E_\xi \subseteq X\) finite so that \(A \subseteq U_\xi\), \(E_\xi \cap V_j \neq \emptyset\) for each \(j\), and \(E_\xi \cap U_\xi = \emptyset\).
    Set \(F_\xi = X \setminus U_\xi\) and notice that \(F_\xi \in \psi^{-1}(\xi)\), \(F_\xi \subseteq X \setminus A\), and since \(E_\xi \cap (X \setminus F_\xi) = \emptyset\), \(E_\xi \subseteq F_\xi\).
    Hence, \(\emptyset \neq E_\xi \cap V_j \subseteq F_\xi \cap V_j\) for each \(j\).
    That is, \(F_\xi \in [ A ; V_0 , \ldots , V_m]\) and so \(\mathscr D \in \mathcal D_{\mathbb F(X,\mathcal A)}^{gp}\).

    Now suppose \(\mathscr D = c.\mathscr U \in \mathcal D_{\mathbb F(X,\mathcal A)}^{gp}\).
    Let \(\varphi:\mathscr D \to |\mathscr D|\) be the partition of \(\mathscr D\) into finite sets as specified in the definition of \(\mathcal D_{\mathbb F(X,F)}^{gp}\).
    Notice that \(|\mathscr U| = |\mathscr D|\) and define \(\psi:\mathscr U \to |\mathscr U|\) by \(\psi(U) = \varphi(X \setminus U)\).
    Notice that this partitions \(\mathscr U\) into finite sets and \(U \in \psi^{-1}(\xi)\) if and only if \(X \setminus U \in \varphi^{-1}(\xi)\).
    Let \(A \in \mathcal A\) and \(V_1 , V_2 , \ldots, V_m \subseteq X\) be open so that \(V_j \cap (X \setminus A) \neq\emptyset\) for each \(j\).
    Notice that the basic open set \([A;V_1, \ldots, V_m]\) is non-empty.
    So we can find \(\xi_0 < |\mathscr U|\) so that for all \(\xi \geq \xi_0\), there are \(F_\xi \in \mathscr \varphi^{-1}(\xi) \cap [A;V_1, \ldots, V_m]\).
    Then \(F_\xi \subseteq X \setminus A\) and \(F_\xi \cap V_j \neq\emptyset\) for each \(j\).
    Immediately, we see that \(A \subseteq X \setminus F_\xi\).
    Set \(U_\xi = X \setminus F_\xi\) and notice that \(U_\xi \in \psi^{-1}(\xi)\).
    Let \(x_{j,\xi} \in F_\xi \cap V_j\) and \(E_\xi = \{ x_{1,\xi} , x_{2,\xi} , \ldots , x_{m,\xi} \}\).
    Clearly, \(E_\xi \cap V_j \neq \emptyset\) for each \(j\) and, since \(E_\xi \subseteq F_\xi\), \(E_\xi \cap U_\xi = \emptyset\).
    That is, \(\mathscr U \in \mathcal O_F^{gp}(X,\mathcal A)\).
\end{proof}

\begin{lemma} \label{lemma:ACoverFellWeakGroup}
    Let \(\mathcal A\) be an ideal of closed subsets of \(X\).
    \(\mathscr U \in \mathcal O^{wgp}_F(X, \mathcal A)\) if and only if \(c.\mathscr U \in \mathcal D^{wgp}_{\mathbb F(X,\mathcal A)}\).
\end{lemma}
\begin{proof}
    Suppose \(\mathscr U \in \mathcal O_F^{wgp}(X,\mathcal A)\).
    Let \(\varphi:\mathscr U \to |\mathscr U|\) be the partition of \(\mathscr U\) into finite sets as specified in the definition of \(\mathcal O_F^{wgp}(X,\mathcal A)\).
    Set \(\mathscr D = c.\mathscr U\) and, noting that \(|\mathscr D| = |\mathscr U|\), define \(\psi:\mathscr D \to |\mathscr U|\) by \(\psi(F) = \varphi(X \setminus F)\).
    Notice that this partitions \(\mathscr D\) into finite sets and \(F \in \psi^{-1}(\xi)\) if and only if \(X \setminus F \in \varphi^{-1}(\xi)\).
    Let \([A;V_1,\ldots,V_m]\) be a non-empty basic open set in \(\mathbb F(X,\mathcal A)\). Then by Proposition \ref{prop:VietorisNonempty}, \(V_j \cap (X \setminus A) \neq \emptyset\) for each \(j\).
    Thus we can find a \(\xi < |\mathscr U|\) and \(E \subseteq X\) finite so that \(A \subseteq \bigcup \varphi^{-1}(\xi)\), \(E \cap V_j \neq \emptyset\) for each \(j\), and \(E \cap \bigcup \varphi^{-1}(\xi) = \emptyset\).
    Set
    \[
    F = X \setminus \left(\bigcup \{G : G \in \phi^{-1}(\xi)\}\right) = \bigcap \{X \setminus G : G \in \phi^{-1}(\xi)\} = \bigcap \{G : G \in \psi^{-1}(\xi)\}.
    \]
   Notice that \(F \subseteq X \setminus A\), and since \(E \cap (X \setminus F) = \emptyset\), \(E \subseteq F\).
    Hence, \(\emptyset \neq E \cap V_j \subseteq F \cap V_j\) for each \(j\).
    That is, \(F \in [ A ; V_0 , \ldots , V_m]\) and so \(\mathscr D \in \mathcal D_{\mathbb F(X,\mathcal A)}^{wgp}\).

    Now suppose \(\mathscr D = c.\mathscr U \in \mathcal D_{\mathbb F(X,\mathcal A)}^{wgp}\).
    Let \(\varphi:\mathscr D \to |\mathscr D|\) be the partition of \(\mathscr D\) into finite sets as specified in the definition of \(\mathcal D_{\mathbb F(X,\mathcal A)}^{wgp}\).
    Note that \(|\mathscr U| = |\mathscr D|\) and define \(\psi:\mathscr U \to |\mathscr D|\) by \(\psi(U) = \varphi(X \setminus U)\).
    Notice that this partitions \(\mathscr U\) into finite sets and \(U \in \psi^{-1}(\xi)\) if and only if \(X \setminus U \in \varphi^{-1}(\xi)\).
    Let \(A \in \mathcal A\) and \(V_1 , V_2 , \ldots, V_m \subseteq X\) be open so that \(V_j \cap (X \setminus A) \neq\emptyset\) for each \(j\).
    Notice that the basic open set \([A;V_1, \ldots, V_m]\) is non-empty.
    So we can find \(\xi < |\mathscr U|\) so that
    \[
    \bigcap \mathscr \varphi^{-1}(\xi) \in [A;V_1, \ldots, V_m].
    \]
    Set \(F = \bigcap \mathscr \varphi^{-1}(\xi)\).
    Then \(F \subseteq X \setminus A\) and \(F \cap V_j \neq\emptyset\) for each \(j\).
    Immediately, we see that \(A \subseteq X \setminus F\).
    Set \(U = X \setminus F\) and notice that \(U = \bigcup \psi^{-1}(\xi)\).
    Let \(x_j \in F \cap V_j\) and \(E = \{ x_1 , x_2 , \ldots , x_n \}\).
    Clearly, \(E \cap V_j \neq \emptyset\) for each \(j\) and, since \(E \subseteq F\), \(E \cap U = \emptyset\).
    That is, \(\mathscr U \in \mathcal O_F^{wgp}(X,\mathcal A)\).
\end{proof}

\begin{lemma} \label{lem:APiNetwork}
    Fix a space \(X\) and an ideal \(\mathcal A \subseteq \mathbb F(X)\).
    Then \(\mathscr Y \in \Pi^F_{\mathcal A}(X)\) if and only if
    \[
    \{ [A; V_1, \ldots, V_n ] : \langle A, V_1, \ldots, V_n \rangle \in \mathscr Y\} \in \mathcal O_{\mathbb F(X,\mathcal A)}.
    \]
\end{lemma}
\begin{proof}
    Suppose \(\mathscr Y \in \Pi_{\mathcal A}^F(X)\).
    Let \(F \in \mathbb F(X)\) be arbitrary.
    Then we can find \(\langle A, V_1, \ldots, V_n \rangle \in \mathscr Y\) and a finite set \(E \subseteq X\) so that \(E \cap V_j \neq\emptyset\) for each \(j\), \(A \subseteq X \setminus F\), and \(E \cap (X \setminus F) = \emptyset\).
    Notice that \(F \subseteq X \setminus A\) and that \(\emptyset \neq E \cap V_j \subseteq F \cap V_j\) for each \(j\).
    That is, \(F \in [A; V_1, V_2, \ldots, V_n]\).

    Now, suppose
    \[
        \{ [A; V_1, \ldots, V_n ] : \langle A, V_1, \ldots, V_n \rangle \in \mathscr Y\} \in \mathcal O_{\mathbb F(X,\mathcal A)}.
    \]
    Let \(U \subseteq X\) be an arbitrary open set.
    Then we can produce \([A ; V_1 , \ldots, V_n]\) where \(\langle A, V_1, \ldots, V_n \rangle \in \mathscr Y\) so that \(X \setminus U \in [A ; V_1 , \ldots, V_n]\).
    Let \(x_j \in V_j \cap (X\setminus U)\) for each \(j\) and let \(E = \{ x_1, \ldots, x_n \}\).
    Observe that \(A \subseteq U\), \(E \cap V_j \neq \emptyset\) for each \(j\), and \(E \cap U = \emptyset\).
    This finishes the proof.
\end{proof}

\begin{corollary} \label{cor:FullTopEquiv}
    For any space \(X\),
    \begin{enumerate}[label=(\roman*)]
        \item \(\mathscr Y\) if a \(\pi_F\)-network if and only if \(\{ [K; V_1, \ldots, V_n ] : \langle K, V_1, \ldots, V_n \rangle \in \mathscr Y\} \in \mathcal O_{\mathbb F(X,F)}\),
        \item If \(\mathscr Y\) is a \(\pi_V\)-network on \(X\), then \(\{ [ \bigcap_{j=1}^n (X \setminus V_j); V_1, \ldots, V_n ] : \langle V_1, \ldots, V_n \rangle \in \mathscr Y\} \in \mathcal O_{\mathbb F(X,V)}\), and
        \item If \(\{ [A; V_1, \ldots, V_n ] : \langle A, V_1, \ldots, V_n \rangle \in \mathscr Y\} \in \mathcal O_{\mathbb F(X, V)}\), then \[
        \{\langle X \setminus A, V_1, \ldots, V_n \rangle : \langle A, V_1, \ldots, V_n \rangle \in \mathscr Y\}
        \]
        is a \(\pi_V\)-network on \(X\).
    \end{enumerate}
\end{corollary}
\begin{proof}
    This follows immediately from Proposition \ref{prop:NetworksEquiv} and Lemma \ref{lem:APiNetwork}.
\end{proof}

\section{Applications of Cover Lemmas} \label{section:Applications}

\subsection{The Upper Topologies}

\begin{theorem} \label{thm:BigOleKoc}
    Fix a topological space \(X\), \(G \in \mathbb F(X)\), ideals \(\mathcal A\) and \(\mathcal B\) consisting of closed sets, \(\alpha\) an ordinal, and a symbol \(\square \in \{1, \text{fin}\}\).
    Then
    \begin{enumerate}[label=(\roman*)]
        \item \label{thm:upperRothberger}
        \(G^\alpha_\square(\mathcal O(X,\mathcal A), \mathcal O(X,\mathcal B)) \equiv G^\alpha_\square(\mathcal D_{\mathbb F(X, \mathcal A^+)}, \mathcal D_{\mathbb F(X, \mathcal B^+)})\),
        \item \label{thm:upperLocalRothberger}
        \(G^\alpha_\square(\mathcal O(X, X \setminus G, \mathcal A), \mathcal O(X, X \setminus G, \mathcal B)) \equiv G^\alpha_\square(\Omega_{\mathbb F(X, \mathcal A^+),G}, \Omega_{\mathbb F(X, \mathcal B^+),G})\),
        \item \label{thm:upperHurewicz}
        \(G^\alpha_\square(\mathcal O(X,\mathcal A), \mathcal O^{gp}(X,\mathcal B)) \equiv G^\alpha_\square(\mathcal D_{\mathbb F(X, \mathcal A^+)}, \mathcal D^{gp}_{\mathbb F(X, \mathcal B^+)})\),
        \item \label{thm:upperWeaklyGroupable}
        \(G^\alpha_\square(\mathcal O(X,\mathcal A), \mathcal O^{wgp}(X,\mathcal B)) \equiv G^\alpha_\square(\mathcal D_{\mathbb F(X, \mathcal A^+)}, \mathcal D^{wgp}_{\mathbb F(X, \mathcal B^+)})\), and
        \item \label{thm:upperPi} if \(\mathcal B \subseteq \mathcal A\),
        \(G^\alpha_\square(\Pi_{\mathcal A}(X), \Pi_{\mathcal B}(X)) \equiv G^\alpha_\square(\mathcal O_{\mathbb F(X,\mathcal A^+)}, \mathcal O_{\mathbb F(X,\mathcal B^+)})\).
    \end{enumerate}
\end{theorem}
\begin{proof}
    Notice that \(\bigcup \mathcal O(X, \mathcal B) = \bigcup \mathcal O(X,\mathcal A) = \mathscr T_X\) and that \(\bigcup \mathcal D_{\mathbb F(X,\mathcal B^+)} = \bigcup \mathcal D_{\mathbb F(X,\mathcal A^+)} = \mathbb F(X)\).
    Define \(\beta_1:\mathscr T_X \to \mathbb F(X)\) by \(\beta_1(U) = X \setminus U\).
    By Lemma \ref{lemma:ACoverUpperFellDense}, \(\beta_1\) satisfies the conditions of Corollary \ref{cor:CorrespondenceEquivalence} for \(\mathcal O(X,\mathcal A)\), \(\mathcal D_{\mathbb F(X,\mathcal A^+)}\), \(\mathcal O(X,\mathcal B)\), and \(\mathcal D_{\mathbb F(X,\mathcal B^+)}\). This proves \ref{thm:upperRothberger}.

    Equivalences \ref{thm:upperLocalRothberger}, \ref{thm:upperHurewicz}, and \ref{thm:upperWeaklyGroupable} follow in a similar way.
    The function \(\beta_1\) is still used and we instead reference Lemmas \ref{lemma:ACoverUpperFellBlade}, \ref{lemma:ACoverUpperFellGroup}, and \ref{lemma:ACoverUpperFellWeakGroup}.

    To prove part \ref{thm:upperPi}, we first note that \(\bigcup \Pi_{\mathcal A}(X) = \mathcal A\), \(\bigcup \Pi_{\mathcal B}(X) = \mathcal B\), \(\bigcup \mathcal O_{\mathbb F(X,\mathcal A^+)} = \mathscr T_{\mathcal A^+}\), and \(\bigcup \mathcal O_{\mathbb F(X,\mathcal B^+)} = \mathscr T_{\mathcal B^+}\).
    The fact that \(\mathscr T_{\mathcal B^+} \subseteq \mathscr T_{\mathcal A^+}\) is a consequence of Lemma \ref{lem:IdealContainmentTopologyContainment}.
    Define \(\beta_2: \mathcal A \to \mathscr T_{\mathcal A^+}\) by
    \[
    \beta_2(A) = (X \setminus A)^+.
    \]
    We first show that \(\beta_2\) is a bijection.
    Suppose \(A_1, A_2 \in \mathcal A\) and \(A_1 \neq A_2\).
    Without loss of generality, let \(x \in A_1 \setminus A_2\).
    Then \(\{x\} \in (X \setminus A_2)^+\) but \(\{x\} \notin (X \setminus A_1)^+\).
    So \(\beta_2(A_1) \neq \beta_2(A_2)\).
    This shows that \(\beta_2\) is injective.
    The fact that \(\beta_2\) is surjective follows from the fact that every open set in \(\mathbb F(X)\) has the form \((X \setminus A)^+\), because \(\mathcal A\) is an ideal of closed sets.
    The other hypotheses of Corollary \ref{cor:CorrespondenceEquivalence} follow from Lemma \ref{lem:UpperTopPiNetwork}.
\end{proof}

\begin{corollary} \label{MegaKocCor}
    Fix a topological space \(X\) and \(G \in \mathbb F(X)\). Then

    \begin{center}
        \begin{longtable}[c]{|r c l|}
            \hline
            \endfirsthead

            \hline
            \endhead

            \hline
            \endfoot

            \hline
            \endlastfoot

            \multicolumn{3}{|c|}{\textbf{Rothberger-Like Games}} \\
            \hline
            \(G_1(\mathcal O(X, \mathbb K(X)), \mathcal O(X, \mathbb K(X)))\) & \(\equiv\) & \(G_1(\mathcal D_{\mathbb F(X, F^+)}, \mathcal D_{\mathbb F(X, F^+)})\) \\
            \hline
            \(G_1(\mathcal O(X, [X]^{<\omega}), \mathcal O(X, [X]^{<\omega})))\) & \(\equiv\) & \(G_1(\mathcal D_{\mathbb F(X, Z^+)}, \mathcal D_{\mathbb F(X, Z^+)})\) \\
            \hline
            \(G_1(\mathcal O(X, \mathbb K(X)), \mathcal O(X, [X]^{<\omega})))\) & \(\equiv\) & \(G_1(\mathcal D_{\mathbb F(X, F^+)}, \mathcal D_{\mathbb F(X, Z^+)})\) \\
            \hline
            \(G_1(\mathcal O(X, X \setminus G, \mathbb K(X)), \mathcal O(X, X \setminus G, \mathbb K(X)))\) & \(\equiv\) & \(G_1(\Omega_{\mathbb F(X, F^+),G}, \Omega_{\mathbb F(X, F^+),G})\) \\
            \hline
           \(G_1(\mathcal O(X, X \setminus G, [X]^{<\omega}), \mathcal O(X, X \setminus G, [X]^{<\omega}))\) & \(\equiv\) & \(G_1(\Omega_{\mathbb F(X, Z^+),G}, \Omega_{\mathbb F(X, Z^+),G})\) \\
            \hline
            \(G_1(\mathcal O(X, X \setminus G, \mathbb K(X)), \mathcal O(X, X \setminus G, [X]^{<\omega}))\) & \(\equiv\) & \(G_1(\Omega_{\mathbb F(X, F^+),G}, \Omega_{\mathbb F(X, Z^+),G})\) \\
            \hline
            \(G_1(\Pi_k(X), \Pi_k(X))\) & \(\equiv\) & \(G_1(\mathcal O_{\mathbb F(X,F^+)}, \mathcal O_{\mathbb F(X, F^+)})\) \\
            \hline
            \(G_1(\Pi_\omega(X), \Pi_\omega(X))\) & \(\equiv\) & \(G_1(\mathcal O_{\mathbb F(X,Z^+)}, \mathcal O_{\mathbb F(X, Z^+)})\) \\
            \hline
            \multicolumn{3}{|c|}{\textbf{Menger-Like Games}} \\
            \hline
            \(G_{\text{fin}}(\mathcal O(X, \mathbb K(X)), \mathcal O(X, \mathbb K(X)))\) & \(\equiv\) & \(G_{\text{fin}}(\mathcal D_{\mathbb F(X, F^+)}, \mathcal D_{\mathbb F(X, F^+)})\) \\
            \hline
            \(G_{\text{fin}}(\mathcal O(X, [X]^{<\omega}), \mathcal O(X, [X]^{<\omega})))\) & \(\equiv\) & \(G_{\text{fin}}(\mathcal D_{\mathbb F(X, Z^+)}, \mathcal D_{\mathbb F(X, Z^+)})\) \\
            \hline
            \(G_{\text{fin}}(\mathcal O(X, \mathbb K(X)), \mathcal O(X, [X]^{<\omega})))\) & \(\equiv\) & \(G_{\text{fin}}(\mathcal D_{\mathbb F(X, F^+)}, \mathcal D_{\mathbb F(X, Z^+)})\) \\
            \hline
            \(G_{\text{fin}}(\mathcal O(X, X \setminus G, \mathbb K(X)), \mathcal O(X, X \setminus G, \mathbb K(X)))\) & \(\equiv\) & \(G_{\text{fin}}(\Omega_{\mathbb F(X, F^+),G}, \Omega_{\mathbb F(X, F^+),G})\) \\
            \hline
           \(G_{\text{fin}}(\mathcal O(X, X \setminus G, [X]^{<\omega}), \mathcal O(X, X \setminus G, [X]^{<\omega}))\) & \(\equiv\) & \(G_{\text{fin}}(\Omega_{\mathbb F(X, Z^+),G}, \Omega_{\mathbb F(X, Z^+),G})\) \\
            \hline
            \(G_{\text{fin}}(\mathcal O(X, X \setminus G, \mathbb K(X)), \mathcal O(X, X \setminus G, [X]^{<\omega}))\) & \(\equiv\) & \(G_{\text{fin}}(\Omega_{\mathbb F(X, F^+),G}, \Omega_{\mathbb F(X, Z^+),G})\) \\
            \hline
            \(G_{\text{fin}}(\Pi_k(X), \Pi_k(X))\) & \(\equiv\) & \(G_{\text{fin}}(\mathcal O_{\mathbb F(X,F^+)}, \mathcal O_{\mathbb F(X, F^+)})\) \\
            \hline
            \(G_{\text{fin}}(\Pi_\omega(X), \Pi_\omega(X))\) & \(\equiv\) & \(G_{\text{fin}}(\mathcal O_{\mathbb F(X,Z^+)}, \mathcal O_{\mathbb F(X, Z^+)})\) \\
            \hline
            \multicolumn{3}{|c|}{\textbf{Hurewicz-Like Games and Weakly Groupable Games}} \\
            \hline
            \(G_1(\mathcal O(X, \mathbb K(X)), \mathcal O^{gp}(X, \mathbb K(X)))\) & \(\equiv\) & \(G_1(\mathcal D_{\mathbb F(X, F^+)}, \mathcal D^{gp}_{\mathbb F(X, F^+)})\) \\
            \hline
            \(G_1(\mathcal O(X, [X]^{<\omega}), \mathcal O^{gp}(X, [X]^{<\omega}))\) & \(\equiv\) & \(G_1(\mathcal D_{\mathbb F(X, Z^+)}, \mathcal D^{gp}_{\mathbb F(X, Z^+)})\) \\
            \hline
            \(G_{\text{fin}}(\mathcal O(X, \mathbb K(X)), \mathcal O^{gp}(X, \mathbb K(X)))\) & \(\equiv\) & \(G_{\text{fin}}(\mathcal D_{\mathbb F(X, F^+)}, \mathcal D^{gp}_{\mathbb F(X, F^+)})\) \\
            \hline
            \(G_{\text{fin}}(\mathcal O(X, [X]^{<\omega}), \mathcal O^{gp}(X, [X]^{<\omega}))\) & \(\equiv\) & \(G_{\text{fin}}(\mathcal D_{\mathbb F(X, Z^+)}, \mathcal D^{gp}_{\mathbb F(X, Z^+)})\) \\
            \hline
            \(G_1(\mathcal O(X, \mathbb K(X)), \mathcal O^{wgp}(X, \mathbb K(X)))\) & \(\equiv\) & \(G_1(\mathcal D_{\mathbb F(X, F^+)}, \mathcal D^{wgp}_{\mathbb F(X, F^+)})\) \\
            \hline
            \(G_1(\mathcal O(X, [X]^{<\omega}), \mathcal O^{wgp}(X, [X]^{<\omega}))\) & \(\equiv\) & \(G_1(\mathcal D_{\mathbb F(X, Z^+)}, \mathcal D^{wgp}_{\mathbb F(X, Z^+)})\) \\
            \hline
            \(G_{\text{fin}}(\mathcal O(X, \mathbb K(X)), \mathcal O^{wgp}(X, \mathbb K(X)))\) & \(\equiv\) & \(G_{\text{fin}}(\mathcal D_{\mathbb F(X, F^+)}, \mathcal D^{wgp}_{\mathbb F(X, F^+)})\) \\
            \hline
            \(G_{\text{fin}}(\mathcal O(X, [X]^{<\omega}), \mathcal O^{wgp}(X, [X]^{<\omega}))\) & \(\equiv\) & \(G_{\text{fin}}(\mathcal D_{\mathbb F(X, Z^+)}, \mathcal D^{wgp}_{\mathbb F(X, Z^+)})\) \\
            \hline
        \end{longtable}
    \end{center}
\end{corollary}

\begin{notn}
	For \(A \subseteq X\), let \(\mathscr N(A)\) be all open sets \(U\) so that \(A \subseteq U\) and let \(\mathscr N_x = \mathscr N(\{x\})\).
	Set \(\mathscr N[X] = \{\mathscr N_x :x \in X\}\), and in general if \(\mathcal A\) is a collection of subsets of \(X\), then \(\mathscr N[\mathcal A] = \{\mathscr N(A) :A \in \mathcal{A}\}\).
	In the case when \(X\) and \(X^\prime\) represent two topologies on the same underlying set, we will use the notation \(\mathscr N_X(A)\) to denote the collection of open sets relative to the topology according to \(X\) that contain \(A\).
\end{notn}

The following corollary ties this discussion of upper-\(\mathcal A\) topologies to selection games involving the space of continuous functions as witnessed by Theorems 25 and 26 of \cite{CCJH2020}.
Recall that \(G_1(\mathscr N[\mathcal A], \neg \mathcal O(X,\mathcal A))\) is a generalized version of the point-open game.

\begin{corollary}
    Let \(\mathcal A\) be an ideal of closed subsets of a space \(X\).
    The games \(G_1(\mathscr N[\mathcal A], \neg \mathcal O(X,\mathcal A))\) and \(G_1(\mathscr T_{\mathbb F(X,\mathcal A^+)}, \neg \mathcal D_{\mathbb F(X,\mathcal A^+)})\) are equivalent.
\end{corollary}
\begin{proof}
    The duality of \(G_1(\mathscr N[\mathcal A], \neg \mathcal O(X,\mathcal A))\) and \(G_1(\mathcal O(X,\mathcal A), \mathcal O(X,\mathcal A))\) is the conclusion of \cite[Cor. 18]{CCJH2020} as a direct application of \cite[Cor. 26]{ClontzDuality}.
    Then note that \(G_1(\mathcal O(X,\mathcal A), \mathcal O(X,\mathcal A))\) and \(G_1(\mathcal D_{\mathbb F(X,\mathcal A^+)}, \mathcal D_{\mathbb F(X,\mathcal A^+)})\) are equivalent by Theorem \ref{thm:BigOleKoc}.
    The duality between \(G_1(\mathcal D_{\mathbb F(X,\mathcal A^+)}, \mathcal D_{\mathbb F(X,\mathcal A^+)})\) and \(G_1(\mathscr T_{\mathbb F(X,\mathcal A^+)}, \neg \mathcal D_{\mathbb F(X,\mathcal A^+)})\) is a direct application of \cite[Prop. 32]{ClontzDuality}.
\end{proof}

\begin{corollary} \label{cor:KFellEquivalent}
    For any space \(X\),
    \begin{itemize}
        \item the games \(G_1(\mathscr N[\mathbb K(X)], \neg \mathcal O(X, \mathbb K(X)))\) and \(G_1(\mathscr T_{\mathbb F(X,F^+)}, \neg \mathcal D_{\mathbb F(X,F^+)})\) are equivalent, and
        \item the games \(G_1(\mathscr N[[X]^{<\omega}], \neg \mathcal O(X, [X]^{<\omega}))\) and \(G_1(\mathscr T_{\mathbb F(X,Z^+)}, \neg \mathcal D_{\mathbb F(X,Z^+)})\) are equivalent.
    \end{itemize}
\end{corollary}

\subsection{The Full Topologies}

\begin{theorem}
    Fix a topological space \(X\), \(G \in \mathbb F(X)\), ideals \(\mathcal A\) and \(\mathcal B\) consisting of closed sets, \(\alpha\) an ordinal, and a symbol \(\square \in \{1, \text{fin}\}\).
    Then
    \begin{enumerate}[label=(\roman*)]
        \item \label{thm:fullRothberger}
        \(G^\alpha_\square(\mathcal O_F(X,\mathcal A), \mathcal O_F(X,\mathcal B)) \equiv G^\alpha_\square(\mathcal D_{\mathbb F(X, \mathcal A)}, \mathcal D_{\mathbb F(X, \mathcal B)})\),
        \item \label{thm:fullLocalRothberger}
        \(G^\alpha_\square(\mathcal O_F(X, X \setminus G, \mathcal A), \mathcal O_F(X, X \setminus G, \mathcal B)) \equiv G^\alpha_\square(\Omega_{\mathbb F(X, \mathcal A),G}, \Omega_{\mathbb F(X, \mathcal B),G})\),
        \item \label{thm:fullHurewicz}
        \(G^\alpha_\square(\mathcal O_F(X,\mathcal A), \mathcal O^{gp}_F(X,\mathcal B)) \equiv G^\alpha_\square(\mathcal D_{\mathbb F(X, \mathcal A)}, \mathcal D^{gp}_{\mathbb F(X, \mathcal B)})\),
        \item \label{thm:fullWeaklyGroupable}
        \(G^\alpha_\square(\mathcal O_F(X,\mathcal A), \mathcal O^{wgp}_F(X,\mathcal B)) \equiv G^\alpha_\square(\mathcal D_{\mathbb F(X, \mathcal A)}, \mathcal D^{wgp}_{\mathbb F(X, \mathcal B)})\),
        \item \label{thm:fullPi1}
        if \(\mathcal B \subseteq \mathcal A\), then \(G^\alpha_\square(\mathcal O_{\mathbb F(X,\mathcal A)}, \mathcal O_{\mathbb F(X,\mathcal B)}) \leq_{\Two} G^\alpha_\square(\Pi^F_{\mathcal A}(X), \Pi^F_{\mathcal B}(X))\), and
        \item \label{thm:fullPi2}
        if \(\mathcal A \subseteq \mathcal B\), then \(G^\alpha_\square(\Pi^F_{\mathcal A}(X), \Pi^F_{\mathcal B}(X)) \leq_{\Two} G^\alpha_\square(\mathcal O_{\mathbb F(X,\mathcal A)}, \mathcal O_{\mathbb F(X,\mathcal B)})\).
    \end{enumerate}
\end{theorem}
\begin{proof}
    Notice that \(\bigcup \mathcal O(X, \mathcal B) = \bigcup \mathcal O(X,\mathcal A) = \mathscr T_X\) and that \(\bigcup \mathcal D_{\mathbb F(X, \mathcal B)} = \bigcup \mathcal D_{\mathbb F(X,\mathcal A)} = \mathbb F(X)\).
    Define \(\beta:\mathscr T_X \to \mathbb F(X)\) by \(\beta_1(U) = X \setminus U\).
    By Lemma \ref{lem:AFDensity}, \(\beta_1\) satisfies the conditions of Corollary \ref{cor:CorrespondenceEquivalence} for \(\mathcal O_F(X,\mathcal A)\), \(\mathcal D_{\mathbb F(X,\mathcal A)}\), \(\mathcal O_F(X,\mathcal B)\), and \(\mathcal D_{\mathbb F(X,\mathcal B)}\). This proves \ref{thm:fullRothberger}.

    The equivalences \ref{thm:fullLocalRothberger}, \ref{thm:fullHurewicz}, and \ref{thm:fullWeaklyGroupable} follow in a similar way using the same \(\beta\) but instead referencing Lemmas \ref{lem:AFcoversBlades}, \ref{lem:AFGrp}, and \ref{lemma:ACoverFellWeakGroup}.

    To prove part \ref{thm:fullPi1}, we first note that
    \begin{itemize}
        \item \(\bigcup \Pi^F_{\mathcal A}(X) = \mathscr X_{\mathcal A}\), \(\bigcup \Pi^F_{\mathcal B}(X) = \mathscr X_{\mathcal B}\),
        \item \(\bigcup \mathcal O_{\mathbb F(X,\mathcal A)} = \mathscr T_{\mathcal A}\), and \(\bigcup \mathcal O_{\mathbb F(X,\mathcal B)} = \mathscr T_{\mathcal B}\)
    \end{itemize}
    Define \(\varphi: \mathscr X_{\mathcal A} \times \alpha \to \mathscr T_{\mathcal A}\) by
    \[
    \varphi(\langle A, V_1, \ldots, V_n \rangle, \xi) = [A; V_1, \ldots, V_n]
    \]
    We check that \(\varphi\) satisfies the hypotheses of Corollary \ref{cor:Translation}.
    Suppose \(\mathscr Y \in \Pi^F_{\mathcal A}(X)\) and \(\xi \in \alpha\).
    Then by Lemma \ref{lem:APiNetwork},
    \[
    \{\varphi(\langle A, V_1, \ldots, V_n \rangle, \xi) : \langle A, V_1, \ldots, V_n \rangle \in \mathscr Y\} = \{[A; V_1, \ldots, V_n] : A \in \mathscr Y\}
    \]
    is an open cover of \(\mathbb F(X,\mathcal A)\).
    Next suppose that \(\mathscr Y_\xi \in \Pi^F_{\mathcal A}(X)\) for each \(\xi < \alpha\), that \(\mathcal G_\xi\) is a finite subset of \(\mathscr Y_\xi\), and that \(\bigcup_{\xi < \alpha} \varphi[\mathcal G_\xi \times \{\xi\}]\) is an open cover of \(\mathbb F(X,\mathcal B)\).
    We need to show that \(\bigcup_{\xi < \alpha}\mathcal G_\xi \in \Pi^F_{\mathcal B}(X)\).
    Write
    \[
        \mathcal G_\xi = \{\langle A_{\xi,k}, V_{1,k}, \ldots, V_{n(k),k} \rangle : k \leq m(\xi)\}
    \]
    Then \(\{[A_{\xi,k}; V_{1,k}, \ldots, V_{n(k),k}] : k \leq m(\xi) \mbox{ and } \xi \in \alpha\}\) is an open cover of \(\mathbb F(X,\mathcal B)\).
    So by Lemma \ref{lem:APiNetwork}, \(\bigcup_{\xi < \alpha}\mathcal G_\xi \in \Pi^F_{\mathcal B}(X)\).
    Thus \(G^\alpha_\square(\mathcal O_{\mathbb F(X,\mathcal A)}, \mathcal O_{\mathbb F(X,\mathcal B)}) \leq_{\Two} G^\alpha_\square(\Pi^F_{\mathcal A}(X), \Pi^F_{\mathcal B}(X))\).

    Finally, we prove \ref{thm:fullPi2}.
    Define \(\Tonexi:\mathcal O_{\mathbb F(X,\mathcal A)} \to \Pi^F_{\mathcal A}(X)\) by
    \[
    \Tonexi(\mathscr U) = \{\langle A, V_1, \ldots, V_n \rangle : (\exists U \in \mathscr U)[[A; V_1, \ldots, V_n] \subseteq U]\}.
    \]
    We first show that \(\Tonexi(\mathscr U) \in \Pi^F_{\mathcal A}(X)\).
    Notice that \(\mathscr U\) can be refined to a cover \(\mathscr V\) of basic open sets and that \(\Tonexi(\mathscr V)\) a subset of \(\Tonexi(\mathscr U)\).
    By Lemma \ref{lem:APiNetwork}, \(\Tonexi(\mathscr V) \in \Pi^F_{\mathcal A}(X)\).
    So \(\Tonexi(\mathscr U)\) is, as well.
    Define \(\Ttwoxi:\mathscr X_{\mathcal A} \times \mathcal O_{\mathbb F(X,\mathcal A)} \to \mathscr T_{\mathcal B}\) as follows.
    If \(W = [A; V_1, \ldots, V_n]\) is a subset of some \(U\) in \(\mathscr U\), choose \(U_{W,\mathscr U} \in \mathscr U\) so that \(W \subseteq U_{W,\mathscr U}\). Then
    \[
    \Ttwoxi(\langle A, V_1,\ldots, V_n \rangle, \mathscr U) = U_{[A; V_1, \ldots, V_n], \mathscr U}
    \]
    if there is a \(U \in \mathscr U\) so that \([A; V_1, \ldots, V_n] \subseteq U\), and set \(\Ttwoxi(\langle A, V_1,\ldots, V_n \rangle, \mathscr U) = \mathbb F(X)\) otherwise.
    Since \(\mathcal A \subseteq \mathcal B\), Lemma \ref{lem:IdealContainmentTopologyContainment} implies that \(\mathscr T_{\mathcal A} \subseteq \mathscr T_{\mathcal B}\).
    So \(\Ttwoxi\) has the correct co-domain.

    Now suppose \(\mathscr U \in \mathcal O_{\mathbb F(X,\mathcal A)}\) and that \(Y \in \Tonexi(\mathscr U)\).
    Say \(Y = \langle A, V_1, \ldots, V_n \rangle\).
    We define \(W = [A; V_1, \ldots, V_n]\).
    Then, by the definition of \(U_{W,\mathscr U}\), \(\Ttwoxi(Y,\mathscr U) = U_{W,\mathscr U} \in \mathscr U\).

    Finally suppose \(\mathscr U_\xi \in \mathcal O_{\mathbb F(X,\mathcal A)}\) and that \(Y_{\xi, 1}, \ldots , Y_{\xi, n(\xi)} \in \Tonexi(\mathscr U_\xi)\) for \(\xi < \alpha\).
    Say
    \[
    Y_{\xi, k} = \langle A_{\xi, k}, V_{\xi, k, 1}, \ldots, V_{\xi, k, m(\xi, k)} \rangle.
    \]
    Suppose \(\{Y_{\xi, k} : \xi < \alpha \mbox{ and } k < n(\xi)\} \in \Pi^F_{\mathcal B}(X)\).
    Set \(W_{\xi, k} = [A_{\xi,k}; V_{\xi, k, 1}, \ldots, V_{\xi, k, m(\xi, k)}]\).
    Notice that
    \[
    \Ttwoxi(Y_{\xi,k}, \mathscr U_\xi) = U_{W_{\xi,k},\mathscr U_\xi}.
    \]
    By Lemma \ref{lem:APiNetwork}, we know that \(\{W_{\xi, k} : \xi < \alpha \mbox{ and } k < n(\xi)\} \in \mathcal O_{\mathbb F(X, \mathcal B)}\).
    Thus as the \(W_{\xi, k} \subseteq U_{W_{\xi,k},\mathscr U_\xi}\), \(\{\Ttwoxi(Y_{\xi,k}, \mathscr U_\xi) : \xi < \alpha \mbox{ and } k < n(\xi)\} \in \mathcal O_{\mathbb F(X, \mathcal B)}\) as well.
    Therefore Corollary \ref{cor:SimplifiedTranslation} applies and \(G^\alpha_\square(\Pi^F_{\mathcal A}(X), \Pi^F_{\mathcal B}(X)) \leq_{\Two} G^\alpha_\square(\mathcal O_{\mathbb F(X,\mathcal A)}, \mathcal O_{\mathbb F(X,\mathcal B)})\).
\end{proof}

The following corollary generalizes the results from Sections 3, 4, and 5 of \cite{Li2016}.

\begin{corollary} \label{MegaLiCor}
    Fix a topological space \(X\) and \(G \in \mathbb F(X)\). Then

    \begin{center}
        \begin{longtable}[c]{|r c l|}
            \hline
            \endfirsthead

            \hline
            \endhead

            \hline
            \endfoot

            \hline
            \endlastfoot

            \hline
            \multicolumn{3}{|c|}{\textbf{Rothberger-Like Games}} \\
            \hline
            \(G_1(\mathcal O_F(X,\mathbb K(X)), \mathcal O_F(X,\mathbb K(X)))\) & \(\equiv\) &  \(G_1(\mathcal D_{\mathbb F(X, F)}, \mathcal D_{\mathbb F(X, F)})\) \\
            \hline
            \(G_1(\mathcal O_F(X,\mathbb F(X)), \mathcal O_F(X,\mathbb F(X)))\) & \(\equiv\) &  \(G_1(\mathcal D_{\mathbb F(X, V)}, \mathcal D_{\mathbb F(X, V)})\) \\
            \hline
            \(G_1(\mathcal O_F(X,\mathbb F(X)), \mathcal O_F(X,\mathbb K(X)))\) & \(\equiv\) &  \(G_1(\mathcal D_{\mathbb F(X, V)}, \mathcal D_{\mathbb F(X, F)})\) \\
            \hline
            \(G_1(\mathcal O_F(X, X \setminus G, \mathbb K(X)), \mathcal O_F(X,X \setminus G, \mathbb K(X)))\) & \(\equiv\) &  \(G_1(\Omega_{\mathbb F(X, F),G}, \Omega_{\mathbb F(X, F),G})\) \\
            \hline
            \(G_1(\mathcal O_F(X, X \setminus G, \mathbb F(X)), \mathcal O_F(X,X \setminus G, \mathbb F(X)))\) & \(\equiv\) &  \(G_1(\Omega_{\mathbb F(X, V),G}, \Omega_{\mathbb F(X, V),G})\) \\
            \hline
            \(G_1(\mathcal O_F(X, X \setminus G, \mathbb F(X)), \mathcal O_F(X,X \setminus G, \mathbb K(X)))\) & \(\equiv\) &  \(G_1(\Omega_{\mathbb F(X, V),G}, \Omega_{\mathbb F(X, F),G})\) \\
            \hline
            \(G_1(\Pi^F_{\mathbb K(X)}(X), \Pi^F_{\mathbb K(X)}(X))\) & \(\equiv\) &  \(G_1(\mathcal O_{\mathbb F(X,F)}, \mathcal O_{\mathbb F(X, F)})\) \\
            \hline
            \(G_1(\Pi^F_{\mathbb F(X)}(X), \Pi^F_{\mathbb F(X)}(X))\) & \(\equiv\) &  \(G_1(\mathcal O_{\mathbb F(X,V)}, \mathcal O_{\mathbb F(X, V)})\) \\
            \hline
            \multicolumn{3}{|c|}{\textbf{Menger-Like Games}} \\
            \hline
            \(G_{\text{fin}}(\mathcal O_F(X,\mathbb K(X)),\mathcal O_F(X,\mathbb K(X)))\) & \(\equiv\) &  \(G_{\text{fin}}(\mathcal D_{\mathbb F(X, F)}, \mathcal D_{\mathbb F(X, F)})\) \\
            \hline
            \(G_{\text{fin}}(\mathcal O_F(X,\mathbb F(X)),\mathcal O_F(X,\mathbb F(X)))\) & \(\equiv\) &  \(G_{\text{fin}}(\mathcal D_{\mathbb F(X, V)}, \mathcal D_{\mathbb F(X, V)})\) \\
            \hline
            \(G_{\text{fin}}(\mathcal O_F(X,\mathbb F(X)),\mathcal O_F(X,\mathbb K(X)))\) & \(\equiv\) &  \(G_{\text{fin}}(\mathcal D_{\mathbb F(X, V)}, \mathcal D_{\mathbb F(X, F)})\) \\
            \hline
            \(G_{\text{fin}}(\mathcal O_F(X, X \setminus G, \mathbb K(X)),\mathcal O_F(X, X \setminus G,\mathbb K(X)))\) & \(\equiv\) &  \(G_{\text{fin}}(\Omega_{\mathbb F(X, F),G}, \Omega_{\mathbb F(X, F),G})\) \\
            \hline
            \(G_{\text{fin}}(\mathcal O_F(X, X \setminus G, \mathbb F(X)),\mathcal O_F(X, X \setminus G,\mathbb K(X)))\) & \(\equiv\) &  \(G_{\text{fin}}(\Omega_{\mathbb F(X, V),G}, \Omega_{\mathbb F(X, V),G})\) \\
            \hline
            \(G_{\text{fin}}(\mathcal O_F(X, X \setminus G, \mathbb F(X)),\mathcal O_F(X, X \setminus G,\mathbb K(X)))\) & \(\equiv\) &  \(G_{\text{fin}}(\Omega_{\mathbb F(X, V),G}, \Omega_{\mathbb F(X, F),G})\) \\
            \hline
            \(G_{\text{fin}}(\Pi^F_{\mathbb K(X)}(X), \Pi^F_{\mathbb K(X)}(X))\) & \(\equiv\) &  \(G_{\text{fin}}(\mathcal O_{\mathbb F(X,F)}, \mathcal O_{\mathbb F(X, F)})\) \\
            \hline
            \(G_{\text{fin}}(\Pi^F_{\mathbb F(X)}(X), \Pi^F_{\mathbb F(X)}(X))\) & \(\equiv\) &  \(G_{\text{fin}}(\mathcal O_{\mathbb F(X,V)}, \mathcal O_{\mathbb F(X, V)})\) \\
            \hline
            \multicolumn{3}{|c|}{\textbf{Hurewicz-Like Games and Weakly Groupable Games}} \\
            \hline
            \(G_1(\mathcal O_F(X, \mathbb K(X)), \mathcal O_F^{gp}(X, \mathbb K(X)))\) & \(\equiv\) &  \(G_1(\mathcal D_{\mathbb F(X, F)}, \mathcal D^{gp}_{\mathbb F(X, F)})\) \\
            \hline
            \(G_1(\mathcal O_F(X, \mathbb F(X)), \mathcal O_F^{gp}(X, \mathbb F(X)))\) & \(\equiv\) &  \(G_1(\mathcal D_{\mathbb F(X, V)}, \mathcal D^{gp}_{\mathbb F(X, V)})\) \\
            \hline
            \(G_{\text{fin}}(\mathcal O_F(X, \mathbb K(X)), \mathcal O_F^{gp}(X, \mathbb K(X)))\) & \(\equiv\) &  \(G_{\text{fin}}(\mathcal D_{\mathbb F(X, F)}, \mathcal D^{gp}_{\mathbb F(X, F)})\) \\
            \hline
            \(G_{\text{fin}}(\mathcal O_F(X, \mathbb F(X)), \mathcal O_F^{gp}(X, \mathbb F(X)))\) & \(\equiv\) &  \(G_{\text{fin}}(\mathcal D_{\mathbb F(X, V)}, \mathcal D^{gp}_{\mathbb F(X, V)})\) \\
            \hline
            \(G_1(\mathcal O_F(X, \mathbb K(X)), \mathcal O_F^{wgp}(X, \mathbb K(X)))\) & \(\equiv\) &  \(G_1(\mathcal D_{\mathbb F(X, F)}, \mathcal D^{wgp}_{\mathbb F(X, F)})\) \\
            \hline
            \(G_1(\mathcal O_F(X, \mathbb F(X)), \mathcal O_F^{wgp}(X, \mathbb F(X)))\) & \(\equiv\) &  \(G_1(\mathcal D_{\mathbb F(X, V)}, \mathcal D^{wgp}_{\mathbb F(X, V)})\) \\
            \hline
            \(G_{\text{fin}}(\mathcal O_F(X, \mathbb K(X)), \mathcal O_F^{wgp}(X, \mathbb K(X)))\) & \(\equiv\) &  \(G_{\text{fin}}(\mathcal D_{\mathbb F(X, F)}, \mathcal D^{wgp}_{\mathbb F(X, F)})\) \\
            \hline
            \(G_{\text{fin}}(\mathcal O_F(X, \mathbb F(X)), \mathcal O_F^{wgp}(X, \mathbb F(X)))\) & \(\equiv\) &  \(G_{\text{fin}}(\mathcal D_{\mathbb F(X, V)}, \mathcal D^{wgp}_{\mathbb F(X, V)})\) \\
            \hline
        \end{longtable}
    \end{center}
\end{corollary}

\subsection{Tightness and Bar-Ilan Selection Principles} \label{subsection:BarIlan}

Both \cite{KocinacEtAl2005} and \cite{Li2016} applied their analyses of the hyperspace topology to the tightness, set-tightness and T-tightness of the hyperspace.
We recall the definitions of tightness and set-tightness now.

\begin{definition}
    Suppose \(X\) is a topological space, \(x \in X\) and \(\kappa\) is a cardinal. Then the \emph{tightness} of \(X\) at \(x\) is bounded above by \(\kappa\), written \(t(X,x) \leq \kappa\), if whenever \(x \in \overline{A} \setminus A\), there is a \(B \subseteq A\) so that \(|B| \leq \kappa\) and \(x \in \overline{B}\).
\end{definition}

\begin{definition}
    Suppose \(X\) is a topological space, \(x \in X\) and \(\kappa\) is a cardinal. Then the \emph{set-tightness} of \(X\) at \(x\) is bounded above by \(\kappa\), written \(t_s(X,x) \leq \kappa\), if whenever \(x \in \overline{A} \setminus A\), there is a \(B \subseteq A\) and a function \(\varphi:B \to \kappa\) so that \(x \in \overline{B}\), but \(x \notin \overline{\varphi^{-1}(\alpha)}\) for any \(\alpha < \kappa\).
\end{definition}

This kind of partition property for a blade inspires the next definition.

\begin{definition}
    Let \(\mathcal S\) be a collection, \(S \in \mathcal S\), and \(\kappa\) be a regular cardinal.
    Then \(S\) is \emph{\(\kappa\)-breakable} if there exists \(\varphi : S \to \kappa\) so that \(\varphi^{-1}(\xi) \notin \mathcal S\) for each \(\xi < \kappa\).
    Let \(\mathcal S^{\kappa\text{br}}\) denote the class of \(\kappa\)-breakable elements of \(\mathcal S\).
\end{definition}

Note that \(t(X,x) \leq t_s(X,x)\) for any \(x \in X\).
Also, set-tightness is essentially tightness, but modified by the \(\kappa\)-breakable definition.
The countable-fan tightness game \(G(\Omega_{X,x}, \Omega_{X,x})\) can describe countable tightness.
Tightness for \(\kappa > \omega\) could be captured by \(G^\kappa(\Omega_{X,x}, \Omega_{X,x})\).
In a similar way, set-tightness is captured by \(G^{|X|}(\Omega_{X,x}, \Omega^{\kappa br}_{X,x})\).

\begin{proposition} \label{prop:ConstantSetTightness}
Suppose \(X\) is a topological space, \(x \in X\), \(\lambda = |X|\), and \(\kappa\) is a cardinal. Then \(t_s(X,x) \leq \kappa\) if and only if \(\One \underset{\text{cnst}}{\not\uparrow} G_1^\lambda(\Omega_{X,x}, \Omega^{\kappa br}_{X,x})\).
\end{proposition}
\begin{proof}
    First suppose that \(t_s(X,x) \leq \kappa\).
    Let \(A\) be the unique play of an arbitrary constant strategy for player One in the game.
    Note that \(x \in \overline{A} \setminus A\).
    Then there is a \(B \subseteq A\) and function \(\varphi:B \to \kappa\) so that \(x \in \overline{B}\), but \(x \notin \overline{\varphi^{-1}(\alpha)}\) for any \(\alpha < \kappa\).
    Enumerate \(B\) as \(\{b_\xi : \xi < |B|\}\).
    Note that \(|B| \leq \lambda\).
    We construct a winning play of the game for player Two as follows.
    In round \(\xi < |B|\), player Two plays \(b_\xi\), and if necessary, in round \(\xi \geq |B|\), player Two plays \(b_0\).
    Then at the end of the game, player Two has created \(B\), which we know is a winning play for Two.

    Now suppose \(\One \underset{\text{cnst}}{\not\uparrow} G_1^\lambda(\Omega_{X,x}, \Omega^{\kappa br}_{X,x})\).
    Let \(A \subseteq X\) be so that \(x \in \overline{A} \setminus A\).
    Then repeatedly playing \(A\) is a constant strategy in the game, so player Two has a winning counter play.
    Suppose Two wins against \(A\) with \(\langle x_\xi : \xi < \lambda \rangle\).
    Set \(B = \{x_\xi : \xi < \lambda\}\).
    Then \(B \subseteq A\) and \(B\) is \(\kappa\)-breakable.
    Thus \(t_s(X,x) \leq \kappa\).
\end{proof}

Since we are able to discuss tightness and set-tightness as the outcome of a constant strategy for a selection game, we can recast these notions as Bar-Ilan selection principles.
Then the translation theorems apply and provide a quick proof of the connection between tightness properties on the hyperspace and Lindel{\"{o}}f-like properties on the ground space.
To complete this analysis, we fist prove equivalence lemmas for \(\Omega^{\kappa br}\) that are similar to those in Section \ref{section:Covers}.

\begin{lemma} \label{lemma:ACoverUpperFellBreakBlade}
    Let \(\mathcal A\) be an ideal of closed subsets of \(X\) and \(G \in \mathbb F(X)\).
    \(\mathscr U \in \mathcal O^{\kappa br}(X, X \setminus G, \mathcal A)\) if and only if \(c.\mathscr U \in \Omega^{\kappa br}_{\mathbb F(X,\mathcal A^+),G}\).
\end{lemma}
\begin{proof}
    First suppose \(\mathscr U\) is a \(\kappa\)-breakable \(\mathcal A\)-cover of \(X \setminus G\) by open sets from \(X\).
    Suppose \(\varphi:\mathscr U \to \kappa\) witnesses that \(\mathscr U\) is \(\kappa\)-breakable.
    Define \(\psi: c.\mathscr U \to \kappa\) by \(\psi(X \setminus U) = \varphi(U)\).
    Note that by Lemma \ref{lemma:ACoverUpperFellBlade}, since \(\varphi^{-1}(\xi)\) is not an \(\mathcal A\)-cover of \(X \setminus G\) for all \(\xi < \kappa\), \(\psi^{-1}(\xi) \notin \Omega_{\mathbb F(X,\mathcal A^+),G}\).
    However, since, \(\mathscr U\) is an \(\mathcal A\)-cover of \(X \setminus G\), \(c.\mathscr U \in \Omega_{\mathbb F(X,\mathcal A^+),G}\).
    Therefore \(c.\mathscr U \in \Omega^{br}_{\mathbb F(X,\mathcal A^+),G}\).

    The other direction is similar.
\end{proof}

\begin{lemma} \label{lemma:ACoverFellBreakBlade}
    Let \(\mathcal A\) be an ideal of closed subsets of \(X\) and \(G \in \mathbb F(X)\).
    \(\mathscr U \in \mathcal O^{\kappa br}_F(X, X \setminus G, \mathcal A)\) if and only if \(c.\mathscr U \in \Omega^{\kappa br}_{\mathbb F(X,\mathcal A),G}\).
\end{lemma}
\begin{proof}
    Use the partition provided by the definition of \(\kappa\)-breakability and appeal to Lemma \ref{lem:AFcoversBlades}.
\end{proof}

\begin{proposition}
    The following equivalences hold:
    \begin{enumerate}[label=(\roman*)]
        \item \label{prop:tightness}
        Tightness in the hyperspace topologies can be captured with
        \[
            t(\mathbb F(X,\mathcal A^+),G) \leq \kappa \iff {\Omega_{\mathbb F(X,\mathcal A^+),G} \choose \Omega^{\kappa}_{\mathbb F(X,\mathcal A^+),G}} \iff {\mathcal O(X,X \setminus G, \mathcal A) \choose \mathcal O^\kappa(X,X \setminus G, \mathcal A)}
        \]
        and
        \[
        t(\mathbb F(X,\mathcal A),G) \leq \kappa \iff {\Omega_{\mathbb F(X,\mathcal A),G} \choose \Omega^{\kappa}_{\mathbb F(X,\mathcal A),G}} \iff {\mathcal O_F(X,X \setminus G, \mathcal A) \choose \mathcal O^{\kappa}_F(X,X \setminus G,\mathcal A)}
        \]
        \item \label{prop:settightness}
        Set-tightness in the hyperspace topologies can be captured with
        \[
        t_s(\mathbb F(X,\mathcal A^+),G) \leq \kappa \iff {\Omega_{\mathbb F(X,\mathcal A^+),G} \choose \Omega^{\kappa\text{br}}_{\mathbb F(X,\mathcal A^+),G}} \iff {\mathcal O(X,X \setminus G, \mathcal A) \choose \mathcal O^{\kappa\text{br}}(X,X \setminus G, \mathcal A)}
        \]
        and
        \[
        t_s(\mathbb F(X,\mathcal A),G) \leq \kappa \iff {\Omega_{\mathbb F(X,\mathcal A),G} \choose \Omega^{\kappa\text{br}}_{\mathbb F(X,\mathcal A),G}} \iff {\mathcal O_F(X,X \setminus G, \mathcal A) \choose \mathcal O_F^{\kappa\text{br}}(X,X \setminus G,\mathcal A)}
        \]
    \end{enumerate}
\end{proposition}
\begin{proof}
    For \ref{prop:tightness}, let \(Y\) be a topological space and \(y \in Y\).
    First, suppose that \(t(Y,y) \leq \kappa\).
    Consider the game \(G^\kappa_1(\Omega_{Y,y},\Omega_{Y,y})\).
    We show that One does not have a constant winning strategy for this game.
    Indeed, suppose \(A \subseteq X\) is any set so that \(y \in \overline{A} \setminus A\), which represents a constant strategy for One.
    By the tightness of \(Y\) at \(y\), we can find \(B = \{ b_\xi : \xi < \kappa \} \subseteq A\) so that \(y \in \overline{B}\).
    Let Two play \(b_\xi\) in the \(\xi^{\text{th}}\) inning and notice that Two wins.

    On the other hand, if One doesn't have a constant winning strategy in \(G^\kappa_1(\Omega_{Y,y},\Omega_{Y,y})\), then we see that any \(A \subseteq Y\) with \(y \in \overline{A}\setminus A\) admits \(B \subseteq A\) with \(|B| \leq \kappa\) so that \(y \in \overline{B}\).
    That is, \(t(Y,y) \leq \kappa\).

    To finish this portion of the proposition, apply Remarks \ref{rmk:ClontzBarIlan} and \ref{rmk:BarIlanThings}.

    To prove \ref{prop:settightness}, apply Proposition \ref{prop:ConstantSetTightness}, Remarks \ref{rmk:ClontzBarIlan} and \ref{rmk:BarIlanThings}, and Lemmas \ref{lemma:ACoverUpperFellBreakBlade} and \ref{lemma:ACoverFellBreakBlade}.
\end{proof}

\section{Further Work}

In Example \ref{example:GroupProper}, we demonstrated a cover which is in \(\mathcal O(\omega_1, [\omega_1]^{<\omega})\), but not in \(\mathcal O^{gp}(\omega_1, [\omega_1]^{<\omega})\).
This cover, however, is not a countable cover.
We do not currently know if there is a example separating \(\mathcal O(X, \mathcal A)\) from \(\mathcal O^{gp}(X, \mathcal A)\) which is countable.
Based on the results of \cite{KocinacScheepers7}, finding such an example for \(\omega\)-covers would be equivalent to finding a space which is Hurewicz but one of its finite powers is not.

We demonstrated equivalences between games on a space \(X\) and games on the Fell/Vietoris hyperspaces of \(X\).
Another common hyperspace construction creates the Pixley-Roy hyperspace (see \cite{Daniels} and \cite{Sakai2012}).
Here we used the natural mapping of open sets to the complements to connect open sets in \(X\) to points in the hyperspace, and this lifted to various game equivalences.
Is there such a connection between \(X\) and the Pixley-Roy hyperspace of \(X\)?

Can Proposition \ref{prop:RefinementEquivalences} be extended to include a broader variety of cover types like groupable covers and weakly groupable covers?

Finally, we did not exhaust the varieties of tightness covered by \cite{KocinacEtAl2005} and \cite{Li2016}.
They connect T-tightness on the Fell/Vietoris hyperspace to a covering property of \(X\).
We could not find a way to translate T-tightness into a more traditional selection principle, and so could not apply the translation result to recover this connection.
Is there a way to characterize T-tightness that makes it fit into our framework?

\providecommand{\bysame}{\leavevmode\hbox to3em{\hrulefill}\thinspace}
\providecommand{\MR}{\relax\ifhmode\unskip\space\fi MR }
\providecommand{\MRhref}[2]{%
  \href{http://www.ams.org/mathscinet-getitem?mr=#1}{#2}
}
\providecommand{\href}[2]{#2}

\end{document}